# THE TAIL OF THE STATIONARY DISTRIBUTION OF A RANDOM COEFFICIENT AR($q$) MODEL[1]

By Claudia Klüppelberg and Serguei Pergamenchtchikov[2]

*Munich University of Technology and Université de Rouen*

We investigate a stationary random coefficient autoregressive process. Using renewal type arguments tailor-made for such processes, we show that the stationary distribution has a power-law tail. When the model is normal, we show that the model is in distribution equivalent to an autoregressive process with ARCH errors. Hence, we obtain the tail behavior of any such model of arbitrary order.

**1. Introduction.** We consider the following random coefficient autoregressive model:

$$(1.1) \qquad y_n = \alpha_{1n} y_{n-1} + \cdots + \alpha_{qn} y_{n-q} + \xi_n, \qquad n \in \mathbb{N},$$

with random variables (r.v.'s) $\alpha_{in} = a_i + \sigma_i \eta_{in}$, where $a_i \in \mathbb{R}$ and $\sigma_i \geq 0$. Set

$$\alpha_n = (\alpha_{1n}, \ldots, \alpha_{qn})', \qquad \eta_n = (\eta_{1n}, \ldots, \eta_{qn})',$$

where throughout the paper all vectors are column vectors and "$'$" denotes transposition. We suppose that the sequences of coefficient vectors $(\eta_n)_{n \in \mathbb{N}}$ and noise variables $(\xi_n)_{n \in \mathbb{N}}$ are independent and both sequences are i.i.d. with

$$(1.2) \qquad \mathbf{E}\xi_1 = \mathbf{E}\eta_{i1} = 0 \quad \text{and} \quad \mathbf{E}\xi_1^2 = \mathbf{E}\eta_{i1}^2 = 1, \qquad i = 1, \ldots, q.$$

We are interested in the existence of a stationary version of the process $(y_n)_{n \in \mathbb{N}}$, represented by a r.v. $y_\infty$ and its properties. In this paper we investigate the tail behavior

$$(1.3) \qquad \mathbf{P}(y_\infty > t) \qquad \text{as } t \to \infty.$$

Received January 2002; revised May 2003.

[1]Supported in part by the German Science Foundation (Deutsche Forschungsgemeinschaft), SFB 386, Statistical Analysis of Discrete Structures.

[2]Supported in part by the RFFI-Grant 00-01-880.

This is, in particular, the first step for an investigation of the extremal behavior of the corresponding stationary process, which we will study in a forthcoming paper. Stationarity of (1.1) is guaranteed by condition (D0) below. To obtain the asymptotic behavior of the tail of $y_\infty$ we embed $(y_n)_{n \in \mathbb{N}}$ into a multivariate setup.

Set $Y_n = (y_n, \ldots, y_{n-q+1})'$. Then the multivariate process $(Y_n)$ can be considered in the much wider context of random recurrence equations of the type

$$(1.4) \qquad Y_n = A_n Y_{n-1} + \zeta_n, \qquad n \in \mathbb{N},$$

where $(A_n, \zeta_n)_{n \in \mathbb{N}}$ is an i.i.d. sequence, the $A_n$ are i.i.d. random $(q \times q)$-matrices and the $\zeta_n$ are i.i.d. $q$-dimensional vectors. Moreover, for every $n$, the vector $Y_{n-1}$ is independent of $(A_n, \zeta_n)$.

Such equations play an important role in many applications as, for example, in queueing; see [4] and in financial time series; see [8]. See also [5] for an interesting review article with a wealth of examples.

If the Markov process defined in (1.4) has a stationary distribution and $Y$ has this stationary distribution, then certain results are known on the tail behavior of $Y$. In the one-dimensional case ($q = 1$), Goldie [10] has derived the tail behavior of $Y$ in a very elegant way by a renewal type argument: the tail decreases like a power-law. For the multivariate model, Kesten [14] and Le Page [17] show—under certain conditions on the matrices $A_n$—that $t^\lambda \mathbf{P}(x'Y > t)$ is asymptotically equivalent to a renewal function, that is,

$$(1.5) \quad t^\lambda \mathbf{P}(x'Y > t) \sim G(x, t) = \mathbf{E}_x \sum_{n=0}^{\infty} g(x_n, t - v_n) \qquad \text{as } t \to \infty,$$

where "$\sim$" means that the quotient of both sides tends to a positive constant. Note that if we set $x' = (1, 0, \ldots, 0)$, then we obtain again (1.3). Here $g(\cdot, \cdot)$ is some continuous function satisfying condition (4.1), $(x_n)_{n \geq 0}$ and $(v_n)_{n \geq 0}$ are stochastic processes, defined in (1.10) and (1.11).

In model (1.1) we have $\zeta_n = (\xi_n, 0, \ldots, 0)'$ and

$$(1.6) \qquad A_n = \begin{pmatrix} \alpha_{1n} & \cdots & \alpha_{qn} \\ I_{q-1} & & 0 \end{pmatrix}, \qquad n \in \mathbb{N},$$

where $I_{q-1}$ denotes the identity matrix of order $q - 1$.

Standard conditions for the existence of a stationary solution to (1.4) are given in [15] (see also [11]) and require that

$$(1.7) \qquad \mathbf{E} \log^+ |A_1| < \infty \quad \text{and} \quad \mathbf{E} \log^+ |\zeta_1| < \infty$$

and that the top Lyapunov exponent

$$(1.8) \qquad \widetilde{\gamma} = \lim_{n \to \infty} n^{-1} \log |A_1 \cdots A_n| < 0.$$

In our case, conditions (1.7) are satisfied. Moreover, we can replace (1.8) by the following simpler condition; see, for example, [20].

**(D0)** The eigenvalues of the matrix

$$(1.9) \qquad\qquad \mathbf{E}\, A_1 \otimes A_1$$

have moduli less than one, where "⊗" denotes the Kronecker product of matrices.

In the context of model (1.1) under the assumption that, for any $n \geq 1$, $\det(A_n) \neq 0$ a.s., the processes $(x_n)_{n\geq 0}$ and $(v_n)_{n\geq 0}$ are defined as

$$(1.10) \quad x_0 = x \in S, \qquad x'_n = \frac{x'_{n-1} A_n}{|x'_{n-1} A_n|} = \frac{x' A_1 \cdots A_n}{|x' A_1 \cdots A_n|}, \qquad n \in \mathbb{N},$$

and

$$
\begin{aligned}
& v_0 = 0, \qquad v_n = \sum_{i=1}^{n} u_i = \log |x' A_1 \cdots A_n|, \\
(1.11) \quad & u_n = \log |x'_{n-1} A_n|, \qquad n \in \mathbb{N}.
\end{aligned}
$$

Here $|\cdot|$ denotes the Euclidean norm in $\mathbb{R}^q$ and $|A|^2 = \operatorname{tr} A A'$ is the corresponding matrix norm; we denote, furthermore, $S = \{z \in \mathbb{R}^q : |z| = 1\}$ and $\overline{x} = x/|x|$ for $x \neq 0$.

Since GARCH models are commonly used as volatility models, modelling the (positive) standard deviation of a financial time series, Kesten's work can be applied to such models; see, for example, [6]. Kesten [14, 15] proved and applied a key renewal theorem to the right-hand side of (1.5) under certain conditions on the function $g$, the Markov chain $(x_n)_{n\geq 0}$ and the stochastic process $(v_n)_{n\geq 0}$; a special case is the random recurrence model (1.4) with $\mathbf{P}(A_n > 0) = 1$, for all $n \in \mathbb{N}$. By completely different, namely, point process methods, Basrak, Davis and Mikosch [1] show that for a stationary model (1.4)—again with positive matrices $A_n$—the stationary distribution has a (multivariate) regularly varying tail. Some special examples have been worked out as ARCH(1) and GARCH(1,1); see [10, 12, 19].

The random coefficient model (1.1), however, does not necessarily satisfy the positivity condition on the matrices $A_n$; see Section 2 for examples. On the other hand, it is a special case within Kesten's very general framework. Consequently, we derived a new key renewal theorem in the spirit of Kesten's results, but tailor-made for Markov chains with compact state space, general matrices $A_n$ and functions $g$ satisfying condition (4.1) (see [16], Theorem 2.1). We apply this theorem to the random coefficient model (1.1).

The paper is organized as follows. Our main results are stated in Section 2. We give weak conditions implying a power-law tail for the stationary distribution of the random coefficient model (1.1). For the Gaussian model (all random coefficients and noise variables are Gaussian) we show that model (1.1) is in distribution equivalent to an autoregressive model with ARCH errors of the same order as the random coefficient model. Since the limit variable of the random recurrence model (1.6) is obtained by iteration, products of random matrices have to be investigated. This is done in Section 3. In Section 4 we check the sufficient conditions and apply the key renewal theorem from [16] to model (1.1). Some auxiliary results are summarized in the Appendix.

**2. Main results.** Our first result concerns stationarity of the multivariate process $(Y_n)_{n \in \mathbb{N}}$ given by (1.4). We need some notions from Markov process theory, which can be found, for example, in [18]. The following result is an immediate consequence of Theorem 3 of [9].

THEOREM 2.1. *Consider model* (1.1) *with* $A_n$ *given by* (1.6), *and* $\zeta_n = (\xi_n, 0, \ldots, 0)'$. *We assume that the independent sequences* $\{\eta_{in}, 1 \leq i \leq q, n \in \mathbb{N}\}$ *and* $(\xi_n)_{n \in \mathbb{N}}$ *are both i.i.d. satisfying* (1.2) *and that* $\xi_1$ *has a positive density on* $\mathbb{R}$. *If* (D0) *holds, then* $Y_n = (y_n, \ldots, y_{n-q+1})'$ *converges in distribution to*

$$(2.1) \qquad Y = \zeta_1 + \sum_{k=2}^{\infty} A_1 \cdots A_{k-1} \zeta_k.$$

*Moreover,* $(Y_n)_{n \in \mathbb{N}}$ *is uniformly geometric ergodic.*

REMARK 2.2. (i) From (2.1) we obtain

$$(2.2) \qquad Y \stackrel{d}{=} A_1 Y_1 + \zeta_1,$$

where $Y_1 = \zeta_2 + \sum_{k=3}^{\infty} A_2 \cdots A_{k-1} \zeta_k \stackrel{d}{=} Y$ and $Y_1$ is independent of $(A_1, \zeta_1)$.

(ii) Since $\mathbf{E}((A_1 \cdots A_n) \otimes (A_1 \cdots A_n)) = (\mathbf{E}(A_1 \otimes A_1))^n$ condition (D0) guarantees that

$$(2.3) \qquad \mathbf{E}|A_1 \cdots A_n|^2 \leq c e^{-\gamma n}$$

for some constants $c, \gamma > 0$. From this follows that the series in (2.1) converges a.s. and the second moment of $Y$ is finite; see Theorem 4 of [9].

We require the following additional conditions for the distributions of the coefficient vectors $(\eta_n)_{n \in \mathbb{N}}$ and the noise variables $(\xi_n)_{n \in \mathbb{N}}$ in model (1.1).

(D1) The r.v.'s $\{\eta_{in}, 1 \leq i \leq q, n \in \mathbb{N}\}$ are i.i.d. with symmetric continuous positive density $\phi(\cdot)$, which is nonincreasing on $\mathbb{R}_+$ and moments of all order exist.

(D2) For some $m \in \mathbb{N}$ we assume that $\mathbf{E}(\alpha_{11} - a_1)^{2m} = \sigma_1^{2m} \mathbf{E}\eta_{11}^{2m} \in (1, \infty)$. In particular, $\sigma_1 > 0$.

(D3) The r.v.'s $(\xi_n)_{n \in \mathbb{N}}$ are i.i.d. and $\mathbf{E}|\xi_1|^m < \infty$ for all $m \geq 2$.

**(D4)** For every real sequence $(c_k)_{k\in\mathbb{N}}$ with $0 < \sum_{k=1}^{\infty} |c_k| < \infty$, the r.v. $\tau = \sum_{k=1}^{\infty} c_k \xi_k$ has a symmetric density, which is nonincreasing on $\mathbb{R}_+$.

Condition (D4) looks rather awkward and complicated to verify. Therefore, we give a simple sufficient condition, which is satisfied by many distributions. The proof is given in Section A1.

PROPOSITION 2.3. *If the r.v. $\xi_1$ has bounded, symmetric density $f$, which is continuously differentiable with bounded derivative $f' \leq 0$ on $[0,\infty)$, then condition (D4) holds.*

The following is our main result.

THEOREM 2.4. *Consider model (1.1), with $A_n$ given by (1.6), and $\zeta_n = (\xi_n, 0, \ldots, 0)'$. We assume that the sequences $\{\eta_{in}, 1 \leq i \leq q, n \in \mathbb{N}\}$ and $\{\xi_n, n \in \mathbb{N}\}$ are independent, that conditions (D0)–(D4) hold and that $a_q^2 + \sigma_q^2 > 0$. Then the distribution of the vector (2.1) satisfies*

$$\lim_{t\to\infty} t^\lambda \mathbf{P}(x'Y > t) = h(x), \qquad x \in S.$$

*The function $h(\cdot)$ is strictly positive and continuous on $S$ and the parameter $\lambda$ is given as the unique positive solution of*

$$(2.4) \qquad \kappa(\lambda) = 1,$$

*where for some probability measure $\nu$ on $S$*

$$(2.5) \qquad \kappa(\lambda) := \lim_{n\to\infty} (\mathbf{E}|A_1 \cdots A_n|^\lambda)^{1/n} = \int_S \mathbf{E}|x'A_1|^\lambda \nu(dx),$$

*and the solution of (2.4) satisfies $\lambda > 2$.*

The following model describes an important special case.

DEFINITION 2.5. *If in model (1.1) all coefficients and the noise are Gaussian; that is, $\eta_{i1} \sim \mathcal{N}(0,1)$ for $i = 1, \ldots, q$ and $\xi_1 \sim \mathcal{N}(0,1)$, we call the model (1.1) a Gaussian linear random coefficient model.*

The proof of the following result is given in Section A2.

PROPOSITION 2.6. *We assume the Gaussian model (1.1) with $\sigma_1 > 0$. This process satisfies conditions (D1)–(D4). Furthermore, under condition (D0), the conditional correlation matrix of $Y$ is given by*

$$(2.6) \qquad R = \mathbf{E}(YY'|A_i, \ i \geq 1) = B + \sum_{k=2}^{\infty} A_1 \cdots A_{k-1} B A_{k-1}' \cdots A_1',$$

*where*

$$B = \begin{pmatrix} 1 & 0 & \cdots & 0 \\ 0 & 0 & \cdot & 0 \\ \cdot & \cdot & \cdot & \cdot \\ 0 & 0 & \cdot & 0 \end{pmatrix}.$$

*Moreover, $R$ is positive definite a.s., that is, the vector $Y$ is conditionally nondegenerate Gaussian and $\mathbf{E}|Y|^2 < \infty$.*

We show that the Gaussian model is in distribution equivalent to an autoregressive model with uncorrelated Gaussian errors, which we specify as an autoregressive process with ARCH errors, an often used class of models for financial time series.

LEMMA 2.7.  *Define for the same $q \in \mathbb{N}$, $a_i \in \mathbb{R}$, $\sigma_i \geq 0$ as in model* (1.1),

$$(2.7) \quad x_n = a_1 x_{n-1} + \cdots + a_q x_{n-q} + \sqrt{1 + \sigma_1^2 x_{n-1}^2 + \cdots + \sigma_q^2 x_{n-q}^2} \, \varepsilon_n, \quad n \in \mathbb{N},$$

*with the same initial values $(x_0, \ldots, x_{-q+1}) = (y_0, \ldots, y_{-q+1})$ as for the process* (1.1). *Furthermore, let $(\varepsilon_n)_{n \in \mathbb{N}}$ be i.i.d. $\mathcal{N}(0,1)$ r.v.'s with initial values $(x_0, \ldots, x_{-q+1})$ independent of the sequence $(\varepsilon_n)_{n \in \mathbb{N}}$. Then the stochastic processes $(x_n)_{n \geq 0}$ and the Gaussian linear random coefficient model* (1.1) *have the same distribution.*

PROOF.  We can rewrite model (1.1) in the form

$$(2.8) \quad y_n = a_1 y_{n-1} + \cdots + a_q y_{n-q} + \sqrt{1 + \sigma_1^2 y_{n-1}^2 + \cdots + \sigma_q^2 y_{n-q}^2} \, \tilde{\varepsilon}_n, \quad n \in \mathbb{N},$$

where

$$\tilde{\varepsilon}_n = \frac{\xi_n + \sigma_1 y_{n-1} \eta_{1n} + \cdots + \sigma_q y_{n-q} \eta_{qn}}{\sqrt{1 + \sigma_1^2 y_{n-1}^2 + \cdots + \sigma_q^2 y_{n-q}^2}}, \qquad n \in \mathbb{N},$$

are i.i.d. $\mathcal{N}(0,1)$. This can be seen by calculating characteristic functions. $\square$

REMARK 2.8.  (i) Since $\det(A_n) = \alpha_{qn} = a_q + \sigma_q \eta_{qn}$, the condition $a_q^2 + \sigma_q^2 > 0$ and condition (D1) guarantee that $\det(A_n) \neq 0$ a.s.

(ii) For $q = 1$, model (2.7) was investigated in [3] by different, purely analytic methods. Stationarity of the model was shown for $a_1^2 + \sigma_1^2 < 1$. Under quite general conditions on the noise variables, defining

$$(2.9) \qquad \kappa(\lambda) = \mathbf{E}|a_1 + \sigma_1 \varepsilon|^\lambda,$$

the equation $\kappa(\cdot) = 1$ has a unique positive solution $\lambda$ and the tail of the stationary r.v. $x_\infty$ satisfies

$$\lim_{t \to \infty} t^\lambda \mathbf{P}(x_\infty > t) = c.$$

Moreover, this also covers infinite variance cases, that is, $\lambda$ can be any positive value.

(iii) Kesten proved a result similar to Theorem 2.4 for the process (1.4) (see [14], Theorem 6) under the following condition: There exists $m > 0$ such that $\mathbf{E}(\lambda_*)^m \geq 1$, where $\lambda_* = \lambda_{\min}(A_1 A_1')$ is the minimal eigenvalue of $A_1 A_1'$. However, for the matrix of the form (1.6) we calculate

$$\lambda_* = \inf_{|z|=1} z' A_1 A_1' z$$

$$= \inf_{|z|=1} \left\{ \sum_{j=1}^{q-1} (\alpha_{j1} z_1 - z_{j+1})^2 + \alpha_{q1}^2 z_1^2 \right\} \leq \sum_{j=2}^{q} z_j^2 = 1 \qquad \text{a.s.},$$

$$\lambda_* \leq \left( 1 + \sum_{j=1}^{q-1} \alpha_{j1}^2 \right)^{-1} \alpha_{q1}^2 \qquad \text{a.s.}$$

In the Gaussian case, when the $\eta_{in}$ are all i.i.d. $\mathcal{N}(0,1)$ with $\sigma_1 > 0$ the second inequality implies $\mathbf{P}(\lambda_* < 1) > 0$. Therefore $\mathbf{E}(\lambda_*)^m < 1$ for any $m > 0$. This means, however, that Kesten's Theorem 6 does not apply to the Gaussian linear random coefficient model.

**3. Products of random matrices.** In this section we investigate the function $\kappa(\lambda)$ as defined in (2.5) for matrices $(A_j)_{j \in \mathbb{N}}$ presented in (1.6) derived from model (1.1). Notice that $A_1 \cdots A_n \stackrel{d}{=} A_n \cdots A_1$ for all $n \in \mathbb{N}$, since the $A_j$ are i.i.d. Furthermore, for any function $f : \mathbb{R}^q \to \mathbb{R}$, we write $f(x') = f(x)$ for all $x \in \mathbb{R}^q$. For the following lemma we adapted the corresponding proof from [17].

LEMMA 3.1. *Assume that conditions* (D1) *and* (D2) *are satisfied and* $a_q^2 + \sigma_q^2 > 0$. *Then there exists some probability measure* $\nu$ *on* $S$ *such that for every* $\lambda > 0$,

$$\kappa(\lambda) := \lim_{n \to \infty} (\mathbf{E}|A_1 \cdots A_n|^\lambda)^{1/n} = \int_S \mathbf{E}|x' A_1|^\lambda \nu(dx) > 0.$$

PROOF. Denote by $\mathbf{B}(S)$ the set of bounded measurable functions and by $\mathbf{C}(S)$ the set of continuous functions on $S$. Define, for $\lambda > 0$,

$$(3.1) \qquad Q_\lambda : \mathbf{B}(S) \to \mathbf{B}(S) \qquad \text{by } Q_\lambda(f)(x) = \mathbf{E}|x' A_1|^\lambda f(\overline{x' A_1})$$

for $x \in S$ and $f \in \mathbf{B}(S)$, where $\overline{v} = v/|v|$ for $v \neq 0$. Notice that, if $f$ is continuous, then also $Q_\lambda(f)$ is continuous, that is, $Q_\lambda : \mathbf{C}(S) \to \mathbf{C}(S)$. Denote by $\mathcal{P}(S)$ the set of probability measures on $S$. Since $S$ is compact in $\mathbb{R}^q$, $\mathcal{P}(S)$ is a compact convex set with respect to the weak topology. Furthermore, for every probability measure $\sigma \in \mathcal{P}(S)$, we define the measure $T_\sigma \in \mathcal{P}(S)$ by

$$(3.2) \qquad T_\sigma(f) = \int_S f(x) T_\sigma(dx) = \frac{\int_S Q_\lambda(f)(x)\sigma(dx)}{\int_S Q_\lambda(e)(x)\sigma(dx)},$$

where $e(x) \equiv 1$, $f \in \mathbf{B}(S)$. The operator $T_\sigma : \mathcal{P}(S) \to \mathcal{P}(S)$ is continuous with respect to the weak topology and, by the Schauder–Tykhonov theorem (see [7], page 450), there exists a fixpoint $\nu \in \mathcal{P}(S)$ such that $T_\nu = \nu$, that is, $T_\nu(f) = \nu(f)$ for all $f \in \mathbf{B}(S)$. This implies that

$$\int_S Q_\lambda(f)(x)\nu(dx) = \kappa(\lambda) \int_S f(x)\nu(dx),$$

where

$$\kappa(\lambda) = \int_S Q_\lambda(e)(x)\nu(dx).$$

Notice that for all $n \in \mathbb{N}$,

$$(3.3) \qquad \int_S Q_\lambda^{(n)}(f)(x)\nu(dx) = \kappa^n(\lambda) \int_S f(x)\nu(dx).$$

Here $Q^{(n)}$ is the $n$th power of the operator $Q$. From (3.1) follows for every $f \in \mathbf{B}(S)$

$$(3.4) \qquad Q_\lambda^{(n)}(f)(x) = \mathbf{E}|x'A_1 \cdots A_n|^\lambda f(\overline{x'A_1 \cdots A_n}), \qquad x \in S.$$

Therefore, by (3.3) $\kappa^n(\lambda) = \int_S Q_\lambda^{(n)}(e)(x)\nu(dx) = \int_S \mathbf{E}|x'A_1 \cdots A_n|^\lambda \nu(dx)$. This implies that $\kappa^n(\lambda) \leq \mathbf{E}|A_1 \cdots A_n|^\lambda$. On the other hand, we have

$$(3.5) \qquad \kappa^n(\lambda) = \mathbf{E}|A_1 \cdots A_n|^\lambda \int_S |x'B_n|^\lambda \nu(dx),$$

where $B_n = A_1 \cdots A_n / |A_1 \cdots A_n|$. We show that

$$(3.6) \qquad c_* = \inf_{|B|=1} \int_S |x'B|^\lambda \nu(dx) > 0.$$

Indeed [taking into account that $\int_S |x'B|^\lambda \nu(dx)$ is a continuous function of $B$], if $c_* = 0$, there exists $B$ with $|B| = 1$ such that $\int_S |x'B|^\lambda \nu(dx) = 0$, which means that $\nu\{x \in S : x'B \neq 0\} = 0$. Set $\mathcal{N} = \{x \in S : x'B = 0\}$ and $g(x) = \chi_{\mathcal{N}^c}$, where $\mathcal{N}^c = S \setminus \mathcal{N}$ and $\chi_A$ denotes the indicator function of a set $A$. If $\mathcal{N} \neq \varnothing$, there exist vectors $b_1 \neq 0, \ldots, b_l \neq 0$ with $1 \leq l \leq q$, such that

$$\mathcal{N} \subset \{x \in \mathbb{R}^q : x'B = 0\} = \{x \in \mathbb{R}^q : x'b_1 = 0, \ldots, x'b_l = 0\}.$$

Furthermore, by (3.3), we obtain, for all $n \in \mathbb{N}$,

$$\int_S Q_\lambda^{(n)}(g)(x)\nu(dx) = \kappa^n(\lambda) \int_S g(x)\nu(dx) = 0.$$

By (3.4) this implies for $n = 2q + 1$

$$\mathbf{E} \int_S |x'A_1 \cdots A_{2q+1}|^\lambda g(\overline{x'A_1 \cdots A_{2q+1}})\nu(dx)$$

$$= \int_{\mathcal{N}} \mathbf{E}|x'A_1 \cdots A_{2q+1}|^\lambda g(\overline{x'A_1 \cdots A_{2q+1}})\nu(dx)$$

$$= 0.$$

Since $\nu(\mathcal{N}) = 1$, there exists some $x \in \mathcal{N}$ such that $\overline{x'A_1 \cdots A_{2q+1}} \in \mathcal{N}$ a.s., that is, for all $1 \leq j \leq l$,

$$\mathbf{P}(x'A_1 \cdots A_{2q+1}b_j = 0) = 1.$$

By Lemma A.12 this is only possible if $b_j = 0$, for all $1 \leq j \leq l$; that is, if $B = 0$. But this contradicts $|B| = 1$. Thus we obtained (3.6). Consequently,

$$\mathbf{E}|A_n \cdots A_1|^\lambda \geq \kappa^n(\lambda)$$

$$= \mathbf{E}|A_n \cdots A_1|^\lambda \int_S |x'B_n|^\lambda \nu(dx)$$

$$\geq c_* \mathbf{E}|A_n \cdots A_1|^\lambda,$$

that is,

$$\kappa(\lambda) \leq (\mathbf{E}|A_n \cdots A_1|^\lambda)^{1/n} \leq \frac{\kappa(\lambda)}{(c_*)^{1/n}}$$

and from this inequality Lemma 3.1 follows by taking the limit as $n \to \infty$. $\square$

LEMMA 3.2. *Assume that conditions* (D0)–(D2) *are satisfied and* $a_q^2 + \sigma_q^2 > 0$. *Then equation* (2.4) *has a unique positive solution.*

PROOF. Denote $\Psi(n) = A_n \cdots A_1 = (\Psi_{ij}(n))$. Then $\Psi_{11}(n) = (\alpha_{1n} - a_1) \times \Psi_{11}(n-1) + \mu_n$, where $\mu_n = a_1\Psi_{11}(n-1) + \alpha_{2n}\Psi_{21}(n-1) + \cdots + \alpha_{qn}\Psi_{q1}(n-1)$ independent of $\eta_{1n}$. By the binomial formula and condition (D1) (which implies that all odd moments of $\eta$ are equal to zero) we have for arbitrary $m \in \mathbb{N}$ with $C_{2m}^j = \binom{2m}{j}$,

$$\mathbf{E}(\Psi_{11}(n))^{2m} = \sum_{j=0}^m C_{2m}^{2j} \mathbf{E}((\alpha_1(n) - a_1)^{2j})\mathbf{E}((\Psi_{11}(n-1))^{2j}\mu_n^{2(m-j)})$$

$$\geq s(m)\mathbf{E}(\Psi_{11}(n-1))^{2m},$$

where by condition (D2) $s(m) = \mathbf{E}(\alpha_{1n} - a_1)^{2m} > 1$ for some $m > 1$. Thus $\mathbf{E}(\Psi_{11}(n))^{2m} \geq s(m)^n$, that is, $\mathbf{E}|\Psi(n)|^{2m} \geq \mathbf{E}(\Psi_{11}(n))^{2m} \geq s(m)^n$, which implies that

$$\kappa(2m) = \lim_{n \to \infty} (\mathbf{E}|\Psi(n)|^{2m})^{1/n} \geq s(m) > 1.$$

We show now that $\log \kappa(\lambda)$ is convex for all $\lambda > 0$ and, hence, continuous on $\mathbb{R}_+$. To see the convexity, set

$$\varsigma_n(\lambda) = \frac{1}{n} \log \mathbf{E}|\Psi(n)|^\lambda, \qquad \lambda > 0,$$

and recall that $\log \kappa(\lambda) = \lim_{n \to \infty} \varsigma_n(\lambda)$. Then for $\alpha \in (0, 1)$ and $\beta = 1 - \alpha$ we obtain by Hölder's inequality, for $\lambda, \mu > 0$,

$$\varsigma_n(\alpha\lambda + \beta\mu) \leq \alpha\varsigma_n(\lambda) + \beta\varsigma_n(\mu).$$

By Remark 2.2(ii) condition (D0) implies (2.3), which ensures that $\kappa(\mu) < 1$ for all $0 < \mu \leq 2$. Therefore equation (2.4) has a unique positive root. $\quad\square$

The proof of the following lemma is a simplification of Step 2 of Theorem 3 of [15] adapted to model (1.1); see also [17], Step 2 of Proposition 1.2.

LEMMA 3.3. *Assume that conditions* (D1) *and* (D2) *are satisfied and* $a_q^2 + \sigma_q^2 > 0$. *For every* $\lambda > 0$ *there exists a continuous function* $h(\cdot) > 0$ *such that for* $Q_\lambda$ *as defined in* (3.2),

$$(3.7) \qquad\qquad Q_\lambda(h)(x) = \kappa(\lambda)h(x), \qquad x \in S.$$

*The function* $h$ *is unique up to a positive constant. Moreover, for* $q = 1$, *it is independent of* $x$.

PROOF. For $q = 1$ we have $S = \{1, -1\}$ and it is easy to deduce that any solution of (3.7) is constant on $S$. For $q \geq 2$ we first recall the notation of the proof of Lemma 3.1, in particular (3.5) and (3.6). Set, for $\lambda > 0$,

$$s_n(x) = \frac{Q_\lambda^{(n)}(e)(x)}{\kappa^n(\lambda)} = \frac{\mathbf{E}|x'A_1 \cdots A_n|^\lambda}{\kappa^n(\lambda)}, \qquad x \in S.$$

Using (3.5) and (3.6), we obtain $\sup_{x \in S} s_n(x) \leq 1/c_*$.

Notice that for any $(q \times q)$-matrix $A$ and $\lambda > 0$, choosing $\lambda_* = \min(\lambda, 1)$,

$$||x'A|^\lambda - |y'A|^\lambda| \leq \max(1, \lambda)|x - y|^{\lambda_*}|A|^\lambda, \qquad x, y \in S,$$

which implies $|s_n(x) - s_n(y)| \leq (\max(1, \lambda)/c_*)|x - y|^{\lambda_*}$. By the principle of Arzéla–Ascoli there exists a sequence $(n_k)_{k \in \mathbb{N}}$ with $n_k \to \infty$ as $k \to \infty$

and a continuous function $h(\cdot)$, such that $h_k(x) := \sum_{j=1}^{n_k} s_j(x)/n_k \to h(x)$ uniformly for $x \in S$ and

$$Q_\lambda(h)(x) = \lim_{k \to \infty} Q_\lambda(h_k)(x) = \lim_{k \to \infty} \frac{1}{n_k} \sum_{j=1}^{n_k} Q_\lambda(s_j)(x)$$

$$= \lim_{k \to \infty} \frac{\kappa(\lambda)}{n_k} \sum_{j=1}^{n_k} s_{j+1}(x) = \kappa(\lambda)h(x).$$

If $h(x) = 0$, for some $x \in S$, then $Q_\lambda^{(n)}(h)(x) = 0$ for all $n \in \mathbb{N}$, that is,

$$\mathbf{E}|x'A_1 \cdots A_n|^\lambda h(x_n) = 0,$$

where $x'_n = \overline{x'A_1 \cdots A_n}$, which means that $h(x_n) = 0$, $\mathbf{P}$-a.s., for all $n \in \mathbb{N}$. From Lemma A.9, where $\pi(\cdot)$ denotes the invariant measure of the Markov chain $(x_n)_{n \geq 0}$, we conclude

$$\mathbf{E}_x h(x_n) = 0 \qquad \forall\, n \in \mathbb{N} \quad \Longrightarrow \quad \lim_{n \to \infty} \mathbf{E}_x h(x_n) = \int_S h(z)\pi(dz) = 0$$

$$\Longrightarrow \quad \lim_{k \to \infty} \int_S h_k(z)\pi(dz) = \int_S h(z)\pi(dz) = 0.$$

But on the other hand,

$$\int_S h_k(z)\pi(dz) = \frac{1}{n_k} \sum_{j=1}^{n_k} \frac{1}{\kappa^j(\lambda)} \int_S Q_\lambda^{(j)}(e)(z)\pi(dz)$$

$$= \frac{1}{n_k} \sum_{j=1}^{n_k} \frac{1}{\kappa^j(\lambda)} \mathbf{E}|A_1 \cdots A_j|^\lambda \int_S \frac{|z'A_1 \cdots A_j|^\lambda}{|A_1 \cdots A_j|^\lambda} \pi(dz)$$

$$\geq c_1 \frac{1}{n_k} \sum_{j=1}^{n_k} \frac{\mathbf{E}|A_1 \cdots A_j|^\lambda}{\kappa^j(\lambda)} \geq c_1,$$

where $c_1 = \inf_{|B|=1} \int_S |z'B|^\lambda \pi(dz)$. Assume that $c_1 = 0$. Then there exists a matrix $B$ with $|B| = 1$, such that $\pi(\mathcal{N}^c \cap S) = 0$ for $\mathcal{N} = \{x \in \mathbb{R}^q : x'B = 0\}$. Denote by $\Lambda(\cdot)$ the Lebesgue measure on $S$, then $\Lambda(\mathcal{N} \cap S) = 0$ because $\mathcal{N}$ is a linear subspace of $\mathbb{R}^q$. By Lemma A.9 $\pi$ is equivalent to $\Lambda$; that is, $\pi(\mathcal{N} \cap S) = 0$. This implies that $\pi(S) = \pi(\mathcal{N}^c \cap S) + \pi(\mathcal{N} \cap S) = 0$, which contradicts $\pi(S) = 1$. Hence, $c_1 > 0$ and $h(x) > 0$ for all $x \in S$.

Now assume that there exists some positive function $g \neq h$ satisfying equation (3.7). Define $\Pi_n = A_1 \cdots A_n$. Then for every $n \in \mathbb{N}$, we have

$$g(x) = \frac{Q_\lambda^{(n)}(g)(x)}{\kappa^n(\lambda)} = \frac{\mathbf{E}|x'\Pi_n|^\lambda g(\overline{x'\Pi_n})}{\kappa^n(\lambda)} = \frac{h(x)}{\kappa^n(\lambda)} \tilde{\mathbf{E}}_x f(x'\Pi_n), \qquad x \in S,$$

where $f(z) = g(\overline{z})/h(\overline{z})$, and for every $n \in \mathbb{N}$,

$$\tilde{\mathbf{E}}_x f(x'\Pi_n) = \frac{1}{h(x)} \mathbf{E}|x'\Pi_n|^\lambda h(\overline{x'\Pi_n}) f(x'\Pi_n), \qquad x \in S,$$

that is, $\tilde{\mathbf{E}}_x$ denotes expectation with respect to the measure defined in (4.7). Since the representation for $g$ holds for all $n$ (therefore for $n = 2q+1$), the function $g$ is continuous by Lemma A.14. Define

$$\rho = \sup_{x \in S} \frac{g(x)}{h(x)} = \frac{g(x_0)}{h(x_0)} \quad \text{and} \quad l(x) = \rho h(x) - g(x), \qquad x \in S.$$

Notice that $l(x) \geq 0$ and $l(x_0) = 0$. Next set

$$L(y) = \frac{l(y)}{h(y)} = \frac{Q_\lambda(l)(y)}{\kappa(\lambda)h(y)} = \cdots = \frac{Q_\lambda^{(n)}(l)(y)}{\kappa^n(\lambda)h(y)} = \frac{Q_\lambda^{(n)}(hL)(y)}{\kappa^n(\lambda)h(y)}, \qquad y \in S.$$

We write

$$\sup_{y \in S} L(y) = L(y_0) = \frac{Q_\lambda^{(n)}(hL)(y_0)}{\kappa^n(\lambda)h(y_0)},$$

equivalently, for $x'_n = \overline{y'_0\Pi_n}$, $\mathbf{E}|y'_0\Pi_n|^\lambda h(x_n)L(x_n) = L(y_0)h(y_0)\kappa^n(\lambda)$. Moreover, (3.7) implies that $\mathbf{E}|y'_0\Pi_n|^\lambda h(x_n) = \kappa^n(\lambda)h(y_0)$ for this sequence $(x_n)_{n \geq 0}$ and therefore $\mathbf{E}|y'_0\Pi_n|^\lambda h(x_n)(L(y_0) - L(x_n)) = 0$. Thus, for all $n \in \mathbb{N}$, $L(x_n) = L(y_0)$ $\mathbf{P}$-a.s. and therefore $\mathbf{E}_{y_0}L(x_n) = \mathbf{E}L(\overline{y'_0\Pi_n}) = L(y_0)$. By Lemma A.9, with $\pi(\cdot)$ the invariant measure of $(x_n)_{n \geq 0}$, we get

$$\int_S L(z)\pi(dz) = \lim_{n \to \infty} \mathbf{E}_{y_0}L(x_n) = L(y_0).$$

Since $L(\cdot)$ is continuous and the measure $\pi(\cdot)$ is equivalent to Lebesgue measure, we have

$$L(y_0) = L(z) = L(x_0) = \frac{l(x_0)}{h(x_0)} = 0, \qquad z \in S.$$

Thus $l(z) = 0$ for all $z \in S$ and Lemma 3.3 follows.   $\square$

**4. Renewal theorem for the associated Markov chain.** The next result is based on the renewal theorem in [16] for the stationary Markov chain $(x_n)_{n \geq 0}$ and the processes $(v_n)_{n \geq 0}$ and $(u_n)_{n \geq 1}$ as defined in (1.10) and (1.11), respectively. Some general properties of $(x_n)_{n \geq 0}$ are summarized in Section A4. Let $g : S \times \mathbb{R} \to \mathbb{R}$ be a continuous bounded function satisfying

$$(4.1) \qquad \sum_{l=-\infty}^{\infty} \sup_{x \in S} \sup_{l \leq t \leq l+1} |g(x,t)| < \infty.$$

The renewal theorem in [16] gives the asymptotic behavior of the renewal function

$$G(x,t) = \mathbf{E}_x \sum_{k=0}^{\infty} g(x_k, t - v_k)$$

under the following conditions:

(C1) For the processes $(x_n)_{n \geq 0}$ and $(u_n)_{n \geq 1}$ define the $\sigma$-algebras

$$\mathcal{F}_0 = \sigma\{x_0\}, \qquad \mathcal{F}_n = \sigma\{x_0, x_1, u_1, \ldots, x_n, u_n\}, \qquad n \in \mathbb{N}.$$

Here the initial value $x_0$ is a r.v., which is independent of $(A_n)_{n \in \mathbb{N}}$. For every bounded measurable function $f : \prod_{i=0}^{\infty} (S \times \mathbb{R}) \to \mathbb{R}$ and for every $\mathcal{F}_n$-measurable r.v. $\varrho$,

(4.2)
$$\mathbf{E}(f(\varrho, x_{n+1}, u_{n+1}, \ldots, x_{n+l}, u_{n+l}, \ldots)|\mathcal{F}_n)$$
$$= \mathbf{E}_{x_n} f(\varrho, x_{n+1}, u_{n+1}, \ldots, x_{n+l}, u_{n+l}, \ldots) =: \Phi(x_n, \varrho),$$

that is, $\Phi(x,a) = \mathbf{E}_x f(a, x_1, u_1, \ldots, x_l, u_l, \ldots)$ for all $x \in S$ and $a \in \mathbb{R}$. Moreover, if for $m \in \mathbb{N}$ the function $f : (S \times \mathbb{R})^m \to \mathbb{R}$ is continuous, then $\Phi(x) = \mathbf{E}_x f(x_1, u_1, \ldots, x_m, u_m)$ is continuous on $S$.

(C2) There exists a probability measure $\pi(\cdot)$ on $S$, which is equivalent to Lebesgue measure such that

$$\|\mathbf{P}_x^{(n)}(\cdot) - \pi(\cdot)\| \to 0, \qquad n \to \infty,$$

for all $x \in S$, where $\|\mu\| = \sup_{|f| \leq 1} \int_S f(y) \mu(dy)$ denotes total variation of any measures $\mu$ on $S$. Moreover, there exists a constant $\beta > 0$ such that for all $x \in S$

$$\lim_{n \to \infty} \frac{v_n}{n} = \beta, \qquad \mathbf{P}_x\text{-a.s.}$$

(C3) There exists some number $m \in \mathbb{N}$ such that for all $\nu \in \mathbb{R}$ and for all $\delta > 0$ there exist $y_{\nu,\delta} \in S$ and $\varepsilon_0 = \varepsilon_0(\nu, \delta) > 0$ such that $\forall\, 0 < \varepsilon < \varepsilon_0$

$$\inf_{x \in B_{\delta,\nu}} \mathbf{P}_x(|x_m - y_{\nu,\delta}| < \varepsilon, |v_m - \nu| < \delta) > 0,$$

where $B_{\delta,\nu} = \{x \in S : |x - y_{\nu,\delta}| < \delta\}$.

(C4) There exists some $l \in \mathbb{N}$ such that the function $\Phi_1(x,t) = \mathbf{E}_x \Phi(x_l, v_l, t)$ satisfies

$$\sup_{|x-y| < \varepsilon} \sup_{t \in \mathbb{R}} |\Phi_1(x,t) - \Phi_1(y,t)| \to 0, \qquad \varepsilon \to 0,$$

for every bounded measurable function $\Phi : S \times \mathbb{R} \times \mathbb{R} \to \mathbb{R}$.

Theorem 4.1 ([16]).   *Assume that conditions* (C1)–(C4) *are satisfied. Then for any function $g$ satisfying* (4.1),

$$\lim_{t\to\infty} G(x,t) = \lim_{t\to\infty} \mathbf{E}_x \sum_{k=0}^{\infty} g(x_k, t - v_k) = \frac{1}{\beta} \int_S \pi(dx) \int_{-\infty}^{\infty} g(x,t)\,dt.$$

We apply this renewal theorem to

$$G(x,t) = \frac{1}{e^t} \int_0^{e^t} u^\lambda \mathbf{P}(x'Y > u)\,du, \qquad x \in S, t \in \mathbb{R},$$

where the vector $Y$ is given by (2.1) and $\lambda$ is the unique positive solution of (2.4).

This definition corresponds to an exponential change of measure, equivalently, to an exponential tilting of the bivariate Markov process $(x_n, v_n)_{n\geq 0}$ as follows. Denote by $\tilde{\mathbf{E}}_x$ the expectation with respect to the probability measure $\tilde{\mathbf{P}}_x$, which is defined by

$$
\begin{aligned}
(4.3) \qquad & \tilde{\mathbf{E}}_x F(x_1, u_1, \ldots, x_n, u_n) \\
& = \frac{1}{h(x)} \mathbf{E}|x'A_1 \cdots A_n|^\lambda h(x_n) F(x_1, u_1, \ldots, x_n, u_n)
\end{aligned}
$$

for each measurable function $F$. Then by Kolmogorov's extension $\tilde{\mathbf{P}}$ and $\tilde{\mathbf{E}}$ are the corresponding quantities [as $\mathbf{P}$ and $\mathbf{E}$ are for $(x_n, v_n)_{n\geq 0}$] of the Markov chain $(\tilde{x}_n, \tilde{v}_n)_{n\geq 0}$ defined by the $n$-step transition densities

$$\tilde{p}_{x,v}^{(n)}(dy, dw) = \frac{e^{\lambda w} h(y)}{e^{\lambda v} h(x)} p_{x,v}^{(n)}(dy, dw),$$

where $p_{x,v}^{(n)}(dy, dw)$ is the $n$-step transition density of the original Markov chain $(x_n, v_n)_{n\geq 0}$.

In order to apply Theorem 4.1 we need to check conditions (C1)–(C4).

However, before we treat the general case for arbitrary dimension $q$, we consider the case $q = 1$ in the next example.

Example 4.2.   Consider model (1.1) for $q = 1$ and $0 < a_1^2 + \sigma_1^2 < 1$, then condition (D0) holds. Define $(x_n)_{n\geq 0}$, $(v_n)_{n\geq 0}$ and $(u_n)_{n\in\mathbb{N}}$ as in (1.10) and (1.11), respectively. Assume that conditions (D1) and (D2) are satisfied. In this case the function $\kappa(\cdot)$ is defined by (2.9), and Lemma 3.2 implies that equation $\kappa(\lambda) = 1$ has a unique positive solution. From Lemma 3.3 we conclude that only constant functions satisfy equation (3.7), and we simply set $h(x) = 1$ in (4.3). This case is special in the sense that $S = \{1, -1\}$, that is, the sphere degenerates to two points, and we define the "Lebesgue measure" on $S$ as any point measure with $\Lambda(1) > 0$ and $\Lambda(-1) > 0$. By the ergodic theorem for finite Markov chains one can directly (without Lemma A.9) conclude that the Markov chain $(x_n)_{n \geq 1}$ [defined in (1.10)] is uniformly geometric ergodic with unique invariant distribution $\pi = \tilde{\pi} = (1/2, 1/2)$ with respect to both measures $\mathbf{P}$ and $\tilde{\mathbf{P}}$, that is, the condition (C2) (with respect to $\tilde{\mathbf{P}}$) holds with $\beta = \mathbf{E}|\alpha_{11}|^{\lambda} \log |\alpha_{11}|$, which is positive (cf. [10], Lemma 2.2).

To show condition (C3) for the measure $\tilde{\mathbf{P}}$, set $m = 1$ and $y_{\nu,\delta} = 1$ for $\nu > 0$ and $\delta > 0$. Therefore, taking into account that by condition (D1) the r.v. $\alpha_{11}$ has a positive density, we obtain the inequality in condition (C3) for any $0 < \varepsilon < 1$.

PROPOSITION 4.3. *Consider model* (1.1) *with* $(x_n)_{n \geq 0}$, $(v_n)_{n \geq 0}$ *and* $(u_n)_{n \in \mathbb{N}}$ *defined in* (1.10) *and* (1.11), *respectively. Assume that conditions* (D0)–(D2) *are satisfied and* $a_q^2 + \sigma_q^2 > 0$. *Then conditions* (C1)–(C4) *hold with respect to the measure* $\tilde{\mathbf{P}}_x$ *generated by the finite-dimensional distributions* (4.3).

PROOF. First recall $\Pi_n = A_1 \cdots A_n$ and $x_n' = \overline{x' \Pi_n} = x' \Pi_n / |x' \Pi_n|$ and $v_n' = \log |x' \Pi_n|$. For every bounded measurable function $\Phi(x_n, v_n, t) = f(x' \Pi_n, t)$, with $f(z, t) = \Phi(\overline{z}, \log |z|, t)$, we have by Lemma A.14 immediately that condition (C4) holds.

Next we check (C1). For $n, l \in \mathbb{N}$ we have

$$x_{n+l}' = \frac{x_n' A_{n+1} \cdots A_{n+l}}{|x_n' A_{n+1} \cdots A_{n+l}|} = h_l(x_n, A_{n+1}, \ldots, A_{n+l})$$

and

$$u_{n+l} = \log |x_{n+l-1}' A_{n+l}| = \log \left| \frac{x_n' A_{n+l-1} \cdots A_{n+1}}{|x_n' A_{n+l-1} \cdots A_{n+1}|} A_{n+l} \right|$$
$$= g_l(x_n, A_{n+1}, \ldots, A_{n+l}).$$

Now for every function $f : \prod_{i=0}^{\infty} (S \times \mathbb{R}) \to \mathbb{R}$ and some $\varrho - \mathcal{F}_n$ measurable r.v. $\varrho$ we calculate

$$f(\varrho, x_{n+1}, u_{n+1}, \ldots, x_{n+l}, u_{n+l}, \ldots)$$
$$= f(\varrho, h_1(x_n, A_{n+1}), g_1(x_n, A_{n+1}), \ldots,$$
$$h_l(x_n, A_{n+1}, \ldots, A_{n+l}), g_l(x_n, A_{n+1}, \ldots, A_{n+l}), \ldots)$$
$$= f_1(\varrho, x_n, A_{n+1}, \ldots, A_{n+l}, \ldots).$$

Therefore, $\mathbf{E}(f(\varrho, x_{n+1}, u_{n+1}, \ldots) | \mathcal{F}_n) = \mathbf{E}(f_1(\varrho, x_n, A_{n+1}, \ldots) | \mathcal{F}_n) = \Phi(x_n, \varrho)$, where [notice that $(\varrho, x_n)$ is independent of $(A_{n+1}, \ldots, A_{n+l}, \ldots)$]

$$\Phi(x, a) = \mathbf{E} f_1(a, x, A_{n+1}, \ldots) = \mathbf{E} f_1(a, x, A_1, \ldots)$$
$$= \mathbf{E} f(a, h_1(x, A_1), g_1(x, A_1), \ldots) = \mathbf{E}_x f(a, x_1, u_1, \ldots).$$

This and (4.3) implies for every $m \in \mathbb{N}$ and every bounded function $f_m : \mathbb{R} \times (S \times \mathbb{R})^m \to \mathbb{R}$,

$$(4.4) \qquad \tilde{\mathbf{E}}_x(f_m(\varrho, x_{n+1}, u_{n+1}, \ldots, x_{n+m}, u_{n+m})|\mathcal{F}_n) = \Phi_m(x_n, \varrho),$$

where $\Phi_m(x, a) = \tilde{\mathbf{E}}_x(f_m(a, x_1, u_1, \ldots, x_m, u_m))$.

Denote by $\mu_x$ the measure on the cylindric $\sigma$-algebra $\mathcal{B}$ in $\prod_{i=0}^{\infty}(S \times \mathbb{R})$ generated by the finite-dimensional distributions of $(x_1, u_1, \ldots, x_k, u_k)$ [defined by (4.3) with initial value $x$] on $\mathcal{B}_k$, where $\mathcal{B}_k$ is the Borel $\sigma$-algebra on $(S \times \mathbb{R})^k$ and $\mathcal{B} = \sigma\{\bigcup_{k=1}^{\infty} \mathcal{B}_k\}$. Let furthermore $\mu_{x|\mathcal{F}_n}$ be the conditional (on $\mathcal{F}_n$) infinite-dimensional distribution of $(x_{n+1}, u_{n+1}, \ldots, x_{n+k}, u_{n+k}, \ldots)$. Equality (4.4) implies that the finite-dimensional distributions of the measure $\mu_{x|\mathcal{F}_n}$ coincide with the finite-dimensional distributions of the measure $\mu_{x_n}$; that is, $\mu_{x|\mathcal{F}_n} \equiv \mu_{x_n}$ on $\mathcal{B}$. This implies (4.2) for the measure defined in (4.3). Furthermore, the definitions of $(x_n)_{n\in\mathbb{N}}$ and $(v_n)_{n\in\mathbb{N}}$ imply that for every continuous $f$ also $\Phi(x) = \tilde{\mathbf{E}}_x f(x_1, v_1, \ldots, x_m, v_m)$ is continuous in $x \in S$. Hence condition (C1) holds.

Next we check condition (C2) for $q \geq 2$. The case $q = 1$ has been treated in Example 4.2. We first show

$$(4.5) \qquad \sup_{x \in S} \tilde{\mathbf{E}}_x(\log|x'A_1|)^2 < \infty.$$

To see this notice that for every $\lambda > 0$,

$$\sup_{x \in \mathbb{R}} \frac{|x|^{\lambda}(\log|x|)^2}{1 + |x|^{\lambda+1}} =: c^* < \infty.$$

Hence for every $x \in S$,

$$\tilde{\mathbf{E}}_x(\log|x'A_1|)^2 = \frac{1}{h(x)}\mathbf{E}|x'A_1|^{\lambda}h(\overline{x'A_1})(\log|x'A_1|)^2$$

$$\leq c^* \frac{h^*}{h_*}(1 + \mathbf{E}|A_1|^{\lambda+1}) < \infty,$$

where $h_* = \inf_{x \in S} h(x)$ and $h^* = \sup_{x \in S} h(x)$. This implies (4.5).

Define

$$f(x) = \frac{1}{h(x)}\mathbf{E}|x'A_1|^{\lambda}\log|x'A_1|h(\overline{x'A_1}) = \tilde{\mathbf{E}}_x \log|x'A_1|,$$

and $m_k = \log|x'_{k-1}A_k| - f(x_{k-1})$, then

$$(4.6) \qquad \frac{v_n}{n} = \frac{1}{n}\sum_{k=1}^{n} f(x_{k-1}) + \frac{1}{n}\sum_{k=1}^{n} m_k := \varsigma_n + \frac{1}{n}\sum_{k=1}^{n} m_k.$$

By the strong law of large numbers for square integrable martingales and (4.5) the last term in (4.6) converges to zero $\tilde{\mathbf{P}}_x$-a.s. for any $x \in S$. By Lemma A.9 $(x_n)_{n\in\mathbb{N}}$ is positive Harris recurrent with respect to the measure $\tilde{\mathbf{P}}_x$ as defined in (4.3). Hence we can apply the ergodic theorem to the first term of the right-hand side of (4.6) (see Theorem 17.0.1, page 411 in [18]). This term then converges to the expectation of $f$ with respect to the invariant measure $\tilde{\pi}$:

$$(4.7) \quad \lim_{n\to\infty} \varsigma_n = \beta = \int_S \tilde{\pi}(dz) \frac{1}{h(z)} \mathbf{E}|z'A_1|^\lambda \log|z'A_1|h(\overline{z'A_1}), \qquad \tilde{\pi}\text{-a.s.}$$

This implies

$$\int_S \tilde{\mathbf{P}}_x\left(\lim_{n\to\infty} \frac{v_n}{n} = \beta\right) \tilde{\pi}(dx) = \int_S \tilde{\mathbf{P}}_x\left(\lim_{n\to\infty} \varsigma_n = \beta\right) \tilde{\pi}(dx) = 1.$$

By Lemma A.9 the measure $\tilde{\pi}$ is equivalent to Lebesgue measure, hence

$$(4.8) \qquad \tilde{\mathbf{P}}_x\left(\lim_{n\to\infty} \frac{v_n}{n} = \beta\right) = 1$$

for $\Lambda$-almost all $x \in S$. From condition (C1) we conclude

$$\tilde{\mathbf{P}}_x\left(\lim_{n\to\infty} \frac{v_n}{n} = \beta\right) = \tilde{\mathbf{E}}_x f(x_l, v_l),$$

where $l = 2q + 1$, and

$$f(x,v) = \tilde{\mathbf{P}}_x\left(\lim_{n\to\infty} \frac{v_n + v}{n} = \beta\right).$$

By condition (C4) the function $\tilde{\mathbf{P}}_x(\lim_{n\to\infty} \frac{v_n}{n} = \beta)$ is continuous on $S$ and therefore (4.8) holds for all $x \in S$.

It remains to show that the constant $\beta$ in (4.7) is positive. By (2.3) there exist $c > 0$ and $\gamma > 0$ such that $\mathbf{E}|\Pi_n|^2 \le c e^{-\gamma n}$. Choose $\delta > 0$ such that $d = \gamma - 2\delta > 0$. Then by Chebyshev's inequality,

$$\mathbf{P}(|x'\Pi_n| \ge e^{-\delta n}) \le e^{2\delta n} \mathbf{E}|x'\Pi_n|^2 \le e^{2\delta n} \mathbf{E}|\Pi_n|^2 \le c e^{-dn}.$$

Moreover, for every $0 < \rho < d/\lambda$ and $x'_n = \overline{x'\Pi_n}$, we have

$$\begin{aligned}
\tilde{\mathbf{P}}_x(|x'\Pi_n| < e^{\rho n}) &= h^{-1}(x)\mathbf{E}|x'\Pi_n|^\lambda h(x_n)\chi_{\{|x'\Pi_n|<e^{\rho n}\}} \\
&\le c^*(e^{-\lambda\delta n} + \mathbf{E}|x'\Pi_n|^\lambda \chi_{\{e^{-\delta n}\le|x'\Pi_n|<e^{\rho n}\}}) \\
&\le c^*(e^{-\lambda\delta n} + e^{\lambda\rho n}\mathbf{P}(|x'\Pi_n| \ge e^{-\delta n})) \\
&\le c^*(e^{-\lambda\delta n} + c e^{-(d-\lambda\rho)n}),
\end{aligned}$$

where $c^* = h^*/h_*$, $h^* = \max h$ and $h_* = \min h$. By the lemma of Borel–Cantelli we conclude that for all $x \in S$,

$$\lim_{n\to\infty} \frac{v_n}{n} \ge \rho > 0 \qquad \tilde{\mathbf{P}}_x\text{-a.s.}$$

This verifies condition (C2).

Finally, we check condition (C3) for $q \geq 2$. The case $q = 1$ has already been treated in Example 4.2. We shall show that for $m = 2q + 1$ and $\forall \nu \in \mathbb{R}$, $\forall \delta > 0$, $\forall y \in S$, $\forall \varepsilon > 0$,

$$(4.9) \qquad \inf_{x \in S} \tilde{\mathbf{P}}_x(|x_m - y| < \varepsilon, |v_m - \nu| < \delta) > 0.$$

Indeed, with $L(z) = z/|z|$, consider

$$\tilde{\mathbf{P}}_x(|x_m - y| < \varepsilon, |v_m - \nu| < \delta) = \tilde{\mathbf{P}}_x(x' \Pi_m \in \Gamma_{y,\varepsilon,\delta}),$$

where $\Gamma_{y,\varepsilon,\delta} = \{z \in \mathbb{R}^q \setminus \{0\} : |L(z) - y| < \varepsilon, \, |\log|z| - \nu| < \delta\}$. For every $y \in S$ and every $\nu \in \mathbb{R}$, this set is a nonempty open set in $\mathbb{R}^q$, because the vector $z_0 = e^\nu y \in \Gamma_{y,\varepsilon,\delta}$ ($\forall \nu \in \mathbb{R}$, $\forall \delta > 0$, $\forall y \in S$, $\forall \varepsilon > 0$). This implies that the Lebesgue measure of $\Gamma_{y,\varepsilon,\delta}$ is positive. By Lemma A.13 we conclude that

$$\inf_{x \in S} \tilde{\mathbf{P}}_x(x' \Pi_m \in \Gamma_{y,\varepsilon,\delta}) > 0.$$

This ensures (4.9), which implies condition (C3). $\square$

Define $\tilde{G}(x,t) = G(x,t)/h(x)$, where $h(\cdot) > 0$ satisfies equation (3.7) with positive $\lambda$ for which $\kappa(\lambda) = 1$. Further, recall that by Remark 2.2 $Y \overset{d}{=} A_1 Y_1 + \zeta_1$, where $Y_1 = \zeta_2 + \sum_{k=3}^\infty A_2 \cdots A_{k-1} \zeta_k$ is independent of $(A_1, \zeta_1)$ and $Y_1 \overset{d}{=} Y$. Therefore,

$$(4.10) \qquad \begin{aligned} \tilde{G}(x,t) &= \frac{1}{h(x)e^t} \int_0^{e^t} u^\lambda \mathbf{P}(x' A_1 Y_1 + x' \zeta_1 > u) \, du \\ &=: \tilde{G}_0(x,t) + g(x,t), \end{aligned}$$

where, setting $\tau_1 = x' A_1 Y_1$ and $\tau_2 = x' \zeta_1$,

$$(4.11) \qquad \begin{aligned} \tilde{G}_0(x,t) &= \frac{1}{h(x)e^t} \int_0^{e^t} u^\lambda \mathbf{P}(\tau_1 > u) \, du, \\ g(x,t) &= \frac{1}{h(x)e^t} \int_0^{e^t} u^\lambda \psi(x,u) \, du \end{aligned}$$

with $\psi(x,u) = \mathbf{P}(\tau_1 + \tau_2 > u) - \mathbf{P}(\tau_1 > u)$.

PROPOSITION 4.4. *Assume that conditions* (D0)–(D2) *are satisfied and* $a_q^2 + \sigma_q^2 > 0$. *Then*

$$(4.12) \qquad \tilde{G}(x,t) = \sum_{n=0}^\infty \tilde{\mathbf{E}}_x g(x_n, t - v_n).$$

PROOF. Lemmas 3.1–3.3 ensure the existence of positive solutions of equations (2.4) and (3.7) which are used in the definition of the measure $\tilde{\mathbf{P}}$ in (4.3). Now consider first $\tilde{G}_0(x, t)$. Mapping $u \mapsto u/|x'A_1|$ and using $x'_1 = x'A_1/|x'A_1|$, we obtain

$$\tilde{G}_0(x, t) = \mathbf{E} \frac{|x'A_1|^\lambda}{h(x)e^{t - \log|x'A_1|}} \int_0^{e^t/|x'A_1|} u^\lambda \mathbf{P}(x'_1 Y > u) \, du$$

$$= \tilde{\mathbf{E}}_x \tilde{G}(x_1, t - \log|x'A_1|).$$

Let $\mathbf{B}(S \times \mathbb{R})$ be a linear space of bounded measurable functions $S \times \mathbb{R} \to \mathbb{R}$. Define the linear operator $\Theta \colon \mathbf{B}(S \times \mathbb{R}) \to \mathbf{B}(S \times \mathbb{R})$ by

$$\Theta(f)(x, t) = \tilde{\mathbf{E}}_x f(x_1, t - v_1),$$

where we have used that $v_1 = u_1 = \log|x'A_1|$. Next, recall that by Proposition 4.3, condition (C1) holds for the measure (4.3). This implies that the $n$th power of the operator $\Theta$ is defined by $\Theta^{(n)}(f)(x, t) = \tilde{\mathbf{E}}_x f(x_n, t - v_n)$. Then equation (4.10) translates into the renewal equation $\tilde{G}(x, t) = \Theta(\tilde{G})(x, t) + g(x, t)$ and we obtain for all $n \in \mathbb{N}$ iteratively,

$$\tilde{G}(x, t) = \Theta^{(n)}(\tilde{G})(x, t) + g(x, t) + \Theta(g)(x, t) + \cdots + \Theta^{(n-1)}(g)(x, t).$$

Moreover, condition (D0) implies $\lim_{n \to \infty} \mathbf{E}|\Pi_n| = 0$, giving

$$\Theta^{(n)}(\tilde{G})(x, t) = \tilde{\mathbf{E}}_x \tilde{G}(x_n, t - v_n)$$

$$= \frac{1}{h(x)e^t} \int_0^{e^t} u^\lambda \mathbf{P}(x'\Pi_n Y > u) \, du \to 0, \qquad n \to \infty.$$

This implies (4.12). □

LEMMA 4.5. *Assume the conditions of Theorem* 2.4. *Then for every* $x \in S$, *there exists*

$$\lim_{t \to \infty} G(x, t) = h(x) \frac{1}{\beta} \int_S \tilde{\pi}(dz) \frac{1}{h(z)} \int_0^\infty u^{\lambda - 1} \psi(z, u) \, du$$

$$= h(x)\gamma^* > 0.$$

*Here* $h(\cdot) > 0$ *satisfies equation* (3.7) *with positive* $\lambda$ *for which* $\kappa(\lambda) = 1$, $\beta > 0$ *is defined in* (4.7) *and* $\tilde{\pi}(\cdot)$ *is the stationary measure of the Markov process* $(x_n)_{n \geq 0}$ *under the distribution* $\tilde{\mathbf{P}}$ *as defined in* (4.3).

PROOF. We apply Theorem 4.1 to (4.12). Conditions (C1)–(C4) hold for $q \geq 1$ by Example 4.2 and Proposition 4.3. It remains to show that the function $g$ given by (4.11) satisfies condition (4.1). By Lemma A.10 follows that $g(x, t) \geq 0$ and therefore

$$|g(x, t)| = g(x, t) \leq \frac{1}{h_*}(g_1^*(x, t) + g_2^*(x, t)),$$

where $h_* = \min_{x \in S} h(x)$ and, with $n(t) = e^{\mu t}$ for some $\mu > 0$,

$$g_1^*(x,t) = \frac{1}{e^t} \int_0^{e^t} u^\lambda \mathbf{P}(\tau_1 > u - n(t)) \, du - \frac{1}{e^t} \int_0^{e^t} u^\lambda \mathbf{P}(\tau_1 > u) \, du,$$

$$g_2^*(x,t) = \frac{e^{\lambda t}}{\lambda + 1} \mathbf{P}(\tau_2 > n(t)).$$

We show that these functions satisfy for sufficiently large $t > 0$ the inequality

(4.13) $$g_i^*(x,t) \le c e^{-c_1 t}$$

for constants $c, c_1 > 0$. First notice that it follows immediately from Lemma 3.2 that $\kappa(\theta) < 1$ for every $1 < \theta < \lambda$. Hence by the defintion of $\kappa(\theta)$ in (2.5), for every $\nu \in (\kappa(\theta), 1)$, there exists some $C = C_\nu > 0$ such that for all $n \in \mathbb{N}$,

$$\mathbf{E}|A_1 \cdots A_n|^\theta \le C \nu^n.$$

From this and Hölder's inequality we obtain, for arbitrary $\rho > 0$,

$$\mathbf{E}|\tau_1|^\theta \le \mathbf{E}|A_1|^\theta \mathbf{E}|Y_1|^\theta$$

$$\le 2^{\theta-1} \mathbf{E}|A_1|^\theta \left( \mathbf{E}|\xi_1|^\theta + \mathbf{E} \left( \sum_{k=3}^\infty |A_2 \cdots A_{k-1}| |\xi_k| \right)^\theta \right)$$

$$\le 2^{\theta-1} \mathbf{E}|A_1|^\theta \left( \mathbf{E}|\xi_1|^\theta + C \mathbf{E}|\xi_1|^\theta \sum_{k=3}^\infty \rho^{-\theta(k-2)} \nu^{k-2} \right.$$

$$\left. \times \left( \sum_{k=3}^\infty \rho^{\theta(k-2)/(\theta-1)} \right)^{\theta-1} \right).$$

Now choose in the last term $\rho = \nu^{1/(2\theta)}$. Then for every $1 < \theta < \lambda$, there exists some $m(\theta) > 0$ such that

(4.14) $$\sup_{x \in S} \mathbf{E}|\tau_1|^\theta = \sup_{x \in S} \mathbf{E}|x' A_1 Y_1|^\theta < m(\theta) < \infty.$$

We study now the function $g_1^*(x,t)$. Indeed, for sufficiently large $t > 0$, we have

$$g_1^*(x,t) \le \frac{1}{e^t} \int_0^{e^t - n(t)} (n(t) + u)^\lambda \mathbf{P}(\tau_1 > u) \, du$$

$$- \frac{1}{e^t} \int_0^{e^t} u^\lambda \mathbf{P}(\tau_1 > u) \, du + \frac{(n(t))^{\lambda+1}}{e^t}$$

$$\le c \frac{(n(t))^{\lambda+1}}{e^t}$$

$$+ \frac{1}{e^t} \int_{n(t)}^{e^t - n(t)} u^\lambda ((1 + n(t)u^{-1})^\lambda - 1) \mathbf{P}(\tau_1 > u) \, du$$

$$\leq c\frac{(n(t))^{\lambda+1}}{e^t} + M^* \frac{n(t)}{e^t} \int_{n(t)}^{e^t - n(t)} u^{\lambda - \theta - 1} \, du \, \mathbf{E}|\tau_1|^{\theta}$$

$$\leq c\frac{(n(t))^{\lambda+1}}{e^t} + M^* \frac{m(\theta) n(t)}{\delta e^{(1-\delta)t}}$$

$$\leq c e^{-(1-\mu(\lambda+1))t} + \frac{M^* m(\theta)}{\delta} e^{-(1-\delta-\mu)t},$$

where

$$M^* = \sup_{0 < x \leq 1} ((1+x)^{\lambda} - 1)/x, \qquad \delta = \lambda - \theta \quad \text{and} \quad c = 2^{\lambda} + 1.$$

To obtain (4.13) for the function $g_1^*(x,t)$, choose the parameters $\delta$ and $\mu$ such that $\delta + \mu < 1$ and $0 < \mu < (1 + \lambda)^{-1}$.

The function $g_2^*(x,t)$ satisfies inequality (4.13), because for every $m > 0$ by condition (D3),

$$\sup_{x \in S} \mathbf{E}|\tau_2|^m = \sup_{x \in S} \mathbf{E}|\langle x \rangle_1 \xi_1|^m \leq \mathbf{E}|\xi_1|^m < \infty,$$

where $\langle x \rangle_1$ denotes the first coordinate of $x \in S$. On the other hand, if $t \to -\infty$, we have immediately from definition (4.11),

$$g(x,t) \leq \frac{1}{h_* e^t} \int_0^{e^t} u^{\lambda} \, du \leq \frac{1}{h_*} e^{\lambda t}$$

and, hence, condition (4.1) holds.

Furthermore, taking into account that $\tilde{\pi}$ is equivalent to Lebesgue measure $\Lambda$ on $S$ (see Lemma A.9), by Theorem 4.1 and Lemma A.10 we conclude

$$\lim_{t \to \infty} \frac{G(x,t)}{h(x)} = \lim_{t \to \infty} \tilde{G}(x,t)$$

$$= \frac{1}{\beta} \int_S \tilde{\pi}(dz) \int_{-\infty}^{+\infty} g(z,s) \, ds$$

$$= \frac{1}{\beta} \int_S \tilde{\pi}(dz) \frac{1}{h(z)} \int_0^{+\infty} u^{\lambda-1} \psi(z,u) \, du$$

$$= \gamma^* > 0. \qquad \square$$

The proof of the following lemma is an immediate consequence of the monotone density theorem in regular variation (see, e.g., [2], Theorem 1.7.2).

LEMMA 4.6. *Assume the conditions of Theorem 2.4. Then for every $x \in S$, there exists*

$$\lim_{t \to \infty} t^{\lambda} \mathbf{P}(x'Y > t) = \gamma^* h(x) > 0,$$

*with $h(\cdot)$ and $\gamma^*$ as in Lemma 4.5.*

EXAMPLE 4.7 (Continuation of Example 4.2).    Lemmas 4.5 and 4.6 imply Theorem 2.4 with the limiting constant

$$\gamma^* = \frac{1}{\beta} \int_0^\infty u^{\lambda-1} \frac{(\psi(1,u) + \psi(-1,u))}{2} \, du.$$

Symmetry of the distribution of $\xi$ implies that $\psi(1,u) = \psi(-1,u)$, hence

$$\lim_{t\to\infty} t^\lambda \mathbf{P}(xY > t) = \frac{1}{\beta} \int_0^\infty u^{\lambda-1} (\mathbf{P}(Y > u) - \mathbf{P}(\alpha_{11} Y_1 > u)) \, du$$

for any $x \in S = \{1, -1\}$.

Note that this special case is already covered by Theorem 2.3 of [10].

## APPENDIX

### A.1. A simple sufficient condition for (D4).

PROOF OF PROPOSITION 2.3.    Let $l = \inf\{k \geq 1 : |c_k| > 0\}$. For $n \geq l$, set $\tau_n = \sum_{k=l}^n c_k \xi_k$. If $|c_k| > 0$, then by the condition of this proposition $c_k \xi_k$ has a symmetric density $p_k(\cdot)$, continuously differentiable with derivative $p'_k(\cdot) \leq 0$ on $[0,\infty)$. Therefore $\tau_l$ has a symmetric density, which is non-increasing on $[0,\infty)$. We proceed by induction. Suppose that $\tau_{n-1}$ has a symmetric density $\varphi_{\tau_{n-1}}(\cdot)$, nonincreasing on $[0,\infty)$. We show that $\tau_n$ has a density with these properties. Indeed, if $c_n = 0$, then $\tau_n = \tau_{n-1}$ and we have the same distribution for $\tau_n$. Consider now the case $|c_n| > 0$. By the properties of $p_n(\cdot)$ and of $\varphi_{\tau_{n-1}}(\cdot)$, we can write the density $\varphi_{\tau_n}(\cdot)$ of $\tau_n$ in the following form:

$$\varphi_{\tau_n}(z) = \int_0^\infty p_n(z+u)\varphi_{\tau_{n-1}}(u) \, du + \int_0^z p_n(z-u)\varphi_{\tau_{n-1}}(u) \, du$$
$$+ \int_z^\infty p_n(u-z)\varphi_{\tau_{n-1}}(u) \, du, \qquad z > 0.$$

Therefore the derivative of this function equals

$$\varphi'_{\tau_n}(z) = \int_z^\infty p'_n(u)(\varphi_{\tau_{n-1}}(u-z) - \varphi_{\tau_{n-1}}(u+z)) \, du$$
$$+ \int_0^z p'_n(u)(\varphi_{\tau_{n-1}}(z-u) - \varphi_{\tau_{n-1}}(u+z)) \, du \leq 0, \qquad z > 0,$$

since $p'_n(\cdot) \leq 0$ and $\varphi_{\tau_{n-1}}(\cdot)$ is nonincreasing on $[0,\infty)$. Therefore we obtained that for all $n \geq l$, the r.v. $\tau_n$ has a symmetric continuously differentiable density, which is nonincreasing on $[0,\infty)$. Moreover, since $\tau = \lim_{n\to\infty} \tau_n$ a.s. and the sequence $(\varphi_{\tau_n}(\cdot))_{n \geq l}$ is uniformly bounded, that is,

$$\sup_{z\in\mathbb{R}, n\geq l} \varphi_{\tau_n}(z) \leq \varphi_{\tau_l}(0) < \infty,$$

we have that for every bounded measurable function $g$ with finite support in $\mathbb{R}$

$$\lim_{n \to \infty} \int_{-\infty}^{\infty} g(z) \varphi_{\tau_n}(z) \, dz = \int_{-\infty}^{\infty} g(z) \varphi_\tau(z) \, dz,$$

where $\varphi_\tau(\cdot)$ is the density of $\tau$. Since $\xi_1$ has a continuous density, $\varphi_\tau$ is also continuous. Therefore, for $0 < a < b$, we have for all $0 < \delta < a$,

$$\int_{b-\delta}^{b+\delta} \varphi_\tau(z) \, dz - \int_{a-\delta}^{a+\delta} \varphi_\tau(z) \, dz$$
$$= \lim_{n \to \infty} \left( \int_{b-\delta}^{b+\delta} \varphi_{\tau_n}(z) \, dz - \int_{a-\delta}^{a+\delta} \varphi_{\tau_n}(z) \, dz \right)$$
$$\leq 0.$$

Since $\varphi_\tau(\cdot)$ is continuous, we conclude

$$\varphi_\tau(b) - \varphi_\tau(a) = \lim_{\delta \to 0} \frac{1}{2\delta} \left( \int_{b-\delta}^{b+\delta} \varphi_\tau(z) \, dz - \int_{a-\delta}^{a+\delta} \varphi_\tau(z) \, dz \right) \leq 0. \qquad \square$$

### A.2. Gaussian linear random coefficient models.

Proof of Proposition 2.6. It is evident that conditions (D1)–(D4) hold for this model with $\sigma_1 > 0$, which implies condition (D2).

To show that the conditional correlation matrix (2.6) is positive definite a.s. take some $x \in \mathbb{R}^q$ such that $x'Rx = 0$. Then for $\Pi_k = A_1 \cdots A_k$, $k \in \mathbb{N}$, and $B$ as defined in (2.6),

$$x'Bx + \sum_{k=1}^{\infty} x'\Pi_k B \Pi_k' x = 0.$$

If we denote by $\langle x \rangle_i$ the $i$th coordinate of $x \in \mathbb{R}^q$, the equality above means that $\langle \Pi_k' x \rangle_1 = 0$ for all $k \in \mathbb{N}$. Set $\theta_k(x) = \langle \Pi_k' x \rangle_1$ for $k \in \mathbb{N}$ and $\theta_0(x) = \langle x \rangle_1$. Because of the special form of the matrices (1.6) one can show by induction that

$$(A.1) \quad \theta_k(x) = \begin{cases} \alpha_{1k} \theta_{k-1}(x) + \cdots + \alpha_{k1} \langle x \rangle_1 + \langle x \rangle_{k+1}, & \text{if } 1 \leq k < q, \\ \alpha_{1k} \theta_{k-1}(x) + \cdots + \alpha_{q(k-q+1)} \theta_{k-q}(x), & \text{if } k \geq q. \end{cases}$$

Consequently, if $\theta_k(x) = 0$ for all $0 \leq k \leq q$, then $\langle x \rangle_1 = \cdots = \langle x \rangle_q = 0$. From this we, conclude that $x'Rx = 0$ implies $x = 0$, which means that $R$ is positive definite a.s. $\quad \square$

**A.3. Auxiliary properties of $\Pi_n = A_1 \cdots A_n$.**  We study the asymptotic properties of $\theta_k(x)$ as defined in (A.1). First recall the classical Anderson inequality; see [13], page 214.

LEMMA A.8 (Anderson's inequality).  *Let $\varsigma$ be a r.v. with symmetric continuous density, which is nonincreasing on $[0, \infty)$. Then for every $c \in \mathbb{R}$ and $a > 0$,*

$$\mathbf{P}(|\varsigma + c| \leq a) \leq \mathbf{P}(|\varsigma| \leq a).$$

LEMMA A.9.  *Assume model (1.1), such that conditions* (D1) *and* (D2) *hold and $a_q^2 + \sigma_q^2 > 0$. Then for every $\mu > 0$ and $k \in \mathbb{N}$,*

$$\text{(A.2)} \qquad \lim_{\delta \to 0} \sup_{|\langle x \rangle_1| > \mu} \mathbf{P}(|\theta_k(x)| < \delta) = 0.$$

*Furthermore, for $k = q$ we have*

$$\text{(A.3)} \quad \lim_{\delta \to 0} \sup_{|x| > \mu} \mathbf{P}(|\theta_q(x)| < \delta) = 0, \qquad \lim_{\delta \to 0} \sup_{x \in S} \tilde{\mathbf{P}}_x(|\theta_q(x)| < \delta) = 0,$$

*where $\tilde{\mathbf{P}}$ is defined in* (4.3).

PROOF.  We show first that for $1 \leq j \leq q$ and for every $\epsilon > 0$ such that $\delta/\epsilon \to 0$ as $\delta \to 0$,

$$\text{(A.4)} \qquad \lim_{\delta \to 0} \sup_{x \in \mathbb{R}^q} \mathbf{P}(|\theta_j(x)| < \delta, |\theta_{j-1}(x)| \geq \epsilon) = 0.$$

Recall that $\theta_0(x) = \langle x \rangle_1$. To prove (A.4) notice first that by (A.1)

$$\theta_j(x) = \eta_{1j} \sigma_1 \theta_{j-1}(x) + m_j(x),$$
$$m_j(x) = a_1 \theta_{j-1}(x) + \alpha_{2(j-1)} \theta_{j-2}(x) + \cdots + \alpha_{j1} \langle x \rangle_1 + \langle x \rangle_{j+1} \chi_{\{j < q\}}.$$

Moreover, condition (D2) implies that $\sigma_1 > 0$ and therefore by Anderson's inequality, [taking into account that $\eta_{1j}$ is independent of $\theta_{j-1}(x)$ and $m_j(x)$] we obtain

$$\begin{aligned}
\mathbf{P}(|\theta_j(x)| &< \delta, |\theta_{j-1}(x)| \geq \epsilon) \\
&= \mathbf{P}(|\eta_{1j} \sigma_1 \theta_{j-1}(x) + m_j(x)| < \delta, |\theta_{j-1}(x)| \geq \epsilon) \\
&\leq \mathbf{P}(|\eta_{1j}| < \delta/(\epsilon \sigma_1)).
\end{aligned}$$

From this and condition (D1) we obtain (A.4). Then (A.2) follows by induction.

Next we show (A.3). Introduce for $\delta > 0$ and $1 \leq j \leq q$ the sets $\Gamma_\delta = \bigcap_{j=1}^q \Gamma_{j,\delta}$, where $\Gamma_{j,\delta} = \{|\theta_j(x)| < \epsilon_j\}$ for $\epsilon_j = \epsilon_j(\delta) = \delta^{j/q}$. Notice that (A.4) implies

$$\lim_{\delta \to 0} \sup_{x \in \mathbb{R}^q} \mathbf{P}(\Gamma_{j,\delta} \cap \Gamma_{j-1,\delta}^c) = 0.$$

Set $\alpha^* = \max_{i+j \leq q} |\alpha_{ij}|$ and define $F_\nu = \{|\alpha_{q1}| \geq \nu\}$, $B_N = \{\alpha^* \leq N\}$. Take for any fixed $\nu > 0$, $N > 0$ the set $\Gamma_\delta \cap F_\nu \cap B_N$. The definition of $\theta_j(x)$ in (A.1) implies that on this set $|x| \to 0$ as $\delta \to 0$. Hence, if—as in (A.3)— $|x| \geq \mu$, there exists $\delta_0 = \delta_0(\mu, \nu, N) > 0$ such that $\Gamma_\delta \cap F_\nu \cap B_N = \varnothing$ for all $\delta \leq \delta_0$. Therefore for this $\delta > 0$ and for $x \in \mathbb{R}^q$ with $|x| > \mu$, we obtain

$$\mathbf{P}(|\theta_q(x)| < \delta)$$

$$\leq \mathbf{P}(\Gamma_\delta) + \sum_{j=2}^{q} \mathbf{P}(\Gamma_{j,\delta} \cap \Gamma_{j-1,\delta}^c)$$

$$\leq \mathbf{P}(|\alpha_q(1)| < \nu) + \mathbf{P}(\alpha^* > N) + \sum_{j=2}^{q} \mathbf{P}(\Gamma_{j,\delta} \cap \Gamma_{j-1,\delta})$$

$$\leq \mathbf{P}(|a_q + \sigma_q \eta_{q1}| < \nu) + \frac{\mathbf{E}\alpha^*}{N} + \sum_{j=2}^{q} \mathbf{P}(\Gamma_{j,\delta} \cap \Gamma_{j-1,\delta}).$$

Notice that the conditions $a_q^2 + \sigma_q^2 > 0$ and (D1) guarantee that the first term in the last line tends to zero as $\nu \to 0$. Hence, we obtain the first limiting equality in (A.3). The second equality follows from the first and the definition (4.3). $\square$

In the following lemma we compute the conditional density of $\Pi'_{2q+1} x$ in $\mathbb{R}^q$ with respect to the random vector $\rho = \rho(x) = \Pi'_q x$.

LEMMA A.10. *Assume that* (D1) *and* (D2) *hold,* $a_q^2 + \sigma_q^2 > 0$ *and* $x \neq 0$. *Then the random vector* $\Pi'_{2q+1} x$ *has conditional* $\mathbf{P}$-density $p_1(z|\rho(x)) = f(z, \rho(x))$ *with respect to* $\rho(x)$. *The function* $f(\cdot, \cdot) : \mathbb{R}^q \times \mathbb{R}^q \to [0, \infty)$ *is given by*

$$(A.5) \qquad f(z, y) = \mathbf{E} \frac{1}{|\det T|} p_0(z'T^{-1}, y),$$

*where*

$$(A.6) \qquad T = \begin{pmatrix} \alpha_{1(q+1)} & \alpha_{2(q+1)} & \cdots & \alpha_{q(q+1)} \\ \vdots & \vdots & \vdots & 0 \\ \alpha_{(q-1)3} & \alpha_{q3} & \cdots & 0 \\ \alpha_{q2} & 0 & \cdots & 0 \end{pmatrix}$$

*and for* $z = (z_1, \ldots, z_q) \in \mathbb{R}^q$, $y = (y_1, \ldots, y_q) \in \mathbb{R}^q$

$$p_0(z, y) = \prod_{j=1}^{q} \varphi_j(z_j | z_{j-1}, \ldots, z_1, y),$$

(A.7)    $\varphi_j(z_j|z_{j-1},\ldots,z_1,y) = \chi_{\{|z_{j-1}|>0\}}\mathbf{E}\frac{1}{\sigma_1|z_{j-1}|}\phi\left(\frac{z_j - m_j(z,y)}{\sigma_1 z_{j-1}}\right),$

$\quad\quad m_1(z,y) = a_1 y_1 + y_2, \quad\quad and\ for\ j>1,$

$\quad\quad m_j(z,y) = a_1 z_{j-1} + \alpha_{2(j-1)}z_{j-2}\cdots+\alpha_{j1}y_1 + y_{j+1}\chi_{\{j<q\}},$

*where $z_0 = y_1$ and the density $\phi$ is defined in condition* (D1).

PROOF. Let $x = (x_1,\ldots,x_q)' \in \mathbb{R}^q$ such that $x_q \neq 0$. We show that the vector $\Pi'_{q+1}x$ has density $f(\cdot,x)$ as defined in (A.5). To this end we show first that $x'\Pi_{q+1} = \theta(x)'T$, where the matrix $T$ is defined in (A.6) and $\theta(x) = (\theta_q(x),\ldots,\theta_1(x))' \in \mathbb{R}^q$. By the definition of $A_j$ in (1.6) we have $\langle x'\Pi_{q+1}\rangle_q = \langle x'\Pi_q A_{q+1}\rangle_q = \alpha_{q(q+1)}\langle x'\Pi_q\rangle_1$ and for $1 \leq j \leq q-1$,

$\langle x'\Pi_{q+1}\rangle_j = \langle x'\Pi_q A_{q+1}\rangle_j = \alpha_{j(q+1)}\langle x'\Pi_q\rangle_1 + \langle x'\Pi_q\rangle_{j+1}$

$\quad = \cdots = \alpha_{j(q+1)}\theta_q(x) + \cdots + \alpha_{(q-1)(j+2)}\theta_{j+1}(x) + \alpha_{q(j+1)}\theta_j(x).$

This gives $x'\Pi_{q+1} = \theta(x)'T$. Next note that $a_q^2 + \sigma_q^2 > 0$ implies

$$|\det T| = \prod_{j=1}^q |\alpha_q(j+1)| = \prod_{j=1}^q |a_q + \sigma_q\eta_q(j+1)| > 0, \quad\quad \mathbf{P}\text{-a.s.}$$

Immediately by (A.1) the vector $\theta(x)$ is measurable with respect to $\sigma\{\alpha_{ik}, 1 \leq i \leq q, 1 \leq k \leq q, i+k \leq q+1\}$. Hence, $T$ is independent of $\theta(x)$. Therefore to prove that the vector $\Pi'_{q+1}x$ has density $f(\cdot,x)$, it suffices to prove that $\theta(x)$ has density $p_0(\cdot,x)$ as in (A.7). Indeed, if $x_1 \neq 0$, then condition (D2) guarantees $\sigma_1^2 > 0$ and $\theta_1(x) = \alpha_{11}x_1 + x_2$ has positive density $\varphi_1(\cdot|x)$ as defined in (A.7). This implies that $\theta_1(x) \neq 0$ a.s., and therefore $\theta_2(x) = \alpha_{12}\theta_1(x) + \alpha_{21}x_1 + x_3$ has conditional density $p_{\theta_2}(z_2|\theta_1(x)) = \varphi_2(z_2|\theta_1(x),x)$, where the function $\varphi_2$ is also defined in (A.7). Similarly, we can show that $p_{\theta_j}(z_j|\theta_{j-1}(x),\ldots,\theta_1(x)) = \varphi_j(z_j|\theta_{j-1}(x),\ldots,\theta_1(x),x)$ for every $2 \leq j \leq q$. Therefore $\theta(x) = (\theta_q(x),\ldots,\theta_1(x))'$ has density (A.7) in $\mathbb{R}^q$ provided $x_1 \neq 0$.

To complete the proof we recall that (A.3) implies $\langle\rho(x)\rangle_1 = \theta_q(x) \neq 0$ a.s. for every vector $x \neq 0$. Therefore, taking into account that the $A_n$ are i.i.d. and $\rho(x)$ independent of $\{A_{q+1},\ldots,A_{2q+1}\}$, we obtain that the conditional [with respect to $\rho(x)$] density of the vector $\Pi'_{2q+1}x = (A_{q+1}\cdots A_{2q+1})'\rho(x)$ equals $f(\cdot,\rho(x))$ a.s. for $x \neq 0$.  □

The following result is an immediate consequence of the definition of $\tilde{\mathbf{P}}$ in (4.3) and Lemma A.12.

COROLLARY A.11. *Under the conditions of Lemma* A.10, *the random vector $\Pi'_{2q+1}x$ has a conditional $\tilde{\mathbf{P}}$-density with respect to $\rho(x)$ given by*

(A.8)    $\tilde{p}_1(z|\rho) = \frac{|z|^\lambda h(\overline{z})}{|\rho|^\lambda h(\overline{\rho})}p_1(z|\rho), \quad\quad z,\rho \in \mathbb{R}^q, z \neq 0, \rho \neq 0,$

*for $p_1(z|x)$ as defined in Lemma* A.10.

LEMMA A.12.   *Assume that conditions* (D1) *and* (D2) *hold and $a_q^2 + \sigma_q^2 > 0$. Then for $b, x \in \mathbb{R}^q$ and $x \neq 0$,*

$$\mathbf{P}(x'\Pi_{2q+1}b = 0) > 0 \implies b = 0.$$

PROOF.   Lemma A.10 implies that

$$\mathbf{P}(x'\Pi_{2q+1}b = 0) = \mathbf{E}\mathbf{P}(x'\Pi_{2q+1}b = 0 | \rho(x))$$

$$= \mathbf{E}\int_{\{z \in \mathbb{R}^q \,:\, z'b = 0\}} p_1(z | \rho(x))\, dz.$$

If this probability is positive, then there exists a vector $\rho \in \mathbb{R}^q$ with $\langle \rho \rangle_1 \neq 0$, such that

$$\int_{\{z \in \mathbb{R}^q \,:\, z'b = 0\}} p_1(z | \rho)\, dz > 0.$$

This is possible if and only if $b = 0$ since the Lebesgue measure of the set $\{z \in \mathbb{R}^q : b'z = 0\}$ equals to zero for all $b \neq 0$.  □

Denote by $\mathrm{mes}(\cdot)$ the Lebesgue measure in $\mathbb{R}^q$.

LEMMA A.13.   *Assume that conditions* (D1) *and* (D2) *hold, $q \geq 2$ and $a_q^2 + \sigma_q^2 > 0$. Then there exists some $\delta_0 > 0$ such that for all $0 < \delta < \delta_0$,*

$$\begin{aligned}
(\text{A.9}) \qquad &\inf_{x \in S} \mathbf{P}(x'\Pi_{2q+1} \in B) \geq p_*(\delta)\mu_\delta(B), \\
&\inf_{x \in S} \tilde{\mathbf{P}}_x(x'\Pi_{2q+1} \in B) \geq \tilde{p}_*(\delta)\tilde{\mu}_\delta(B),
\end{aligned}$$

*for every measurable set $B \subseteq \mathbb{R}^q$. Here $p_*(\delta), \tilde{p}_*(\delta) > 0$ and*

$$\begin{aligned}
(\text{A.10}) \qquad &\mu_\delta(B) = \mathbf{E}\int_{\Omega_\delta} \chi_B(z'T)\, dz, \\
&\tilde{\mu}_\delta(B) = \mathbf{E}\int_{\Omega_\delta} |z'T|^\lambda \chi_B(z'T)\, dz, \\
&\Omega_\delta = \{y = (y_1, \ldots, y_q)' \in \mathbb{R}^q : \delta \leq |y_j| \leq \delta^{-1}, j = 1, \ldots, q\},
\end{aligned}$$

*and the matrix $T$ is defined in* (A.6). *Moreover, if $\mathrm{mes}(B) > 0$, then there exists some $\delta_0 > 0$ such that $\mu_\delta(B) > 0$ and $\tilde{\mu}_\delta(B) > 0$ for all $0 < \delta < \delta_0$.*

PROOF.   From Lemma A.10 we know that for a some $0 < \delta < 1$,

$$\mathbf{P}(x'\Pi_{2q+1} \in B) = \mathbf{E}\mathbf{P}(x'\Pi_{2q+1} \in B | \rho(x)) \geq \mathbf{E}\chi_{\{\rho(x) \in K_\delta\}} I_B(\rho(x)),$$

where $K_\delta = \{y = (y_1, \ldots, y_q)' \in \mathbb{R}^q : \delta \leq |y_1| \text{ and } |y| \leq \delta^{-1}\}$ and

$$I_B(\rho) = \int_{\mathbb{R}^q} \chi_B(z) p_1(z|\rho) \, dz$$

$$= \mathbf{E} \int_{\mathbb{R}^q} \chi_B(z'T) p_0(z, \rho) \, dz$$

$$\geq \mathbf{E} \int_{\Omega_\delta} \chi_B(z'T) p_0(z, \rho) \, dz.$$

Next we show for $K_\delta^c = \mathbb{R}^q \setminus K_\delta$,

$$(A.11) \qquad \begin{aligned} &\lim_{\delta \to 0} \sup_{x \in S} \mathbf{P}(\rho(x) \in K_\delta^c) = 0, \\ &\lim_{\delta \to 0} \sup_{x \in S} \tilde{\mathbf{P}}_x(\rho(x) \in K_\delta^c) = 0. \end{aligned}$$

Indeed, we have

$$\mathbf{P}(\rho(x) \in K_\delta^c) \leq \mathbf{P}(|\langle \rho(x) \rangle_1| < \delta) + \mathbf{P}(|\rho(x)| > \delta^{-1})$$

$$\leq \sup_{x \in S} \mathbf{P}(|\theta_q(x)| < \delta) + \delta(\mathbf{E}|A_1|)^q.$$

(A.3) gives the limits in (A.11).

Notice that (A.7) implies that $M_*(\delta) = \inf_{z \in \Omega_\delta, x \in K_\delta} p_0(z, x) > 0$ for every $\delta > 0$, which yields $\mathbf{P}(x'\Pi_{2q+1} \in B) \geq M_*(\delta)\mathbf{P}(\rho(x) \in K_\delta)\mu_\delta(B)$. From this and (A.11) we obtain the first inequality in (A.9). Similarly, we prove the second.

Let $B$ be a measurable set in $\mathbb{R}^q$. By the monotone convergence theorem we have

$$\lim_{\delta \to 0} \mu_\delta(B) = \text{mes}(B)\mathbf{E}|\det T|^{-1},$$

$$\lim_{\delta \to 0} \tilde{\mu}_\delta(B) = \int_{\mathbb{R}^q} |z|^\lambda \chi_B(z) \, dz \, \mathbf{E}|\det T|^{-1}.$$

Since $|\det T| < \infty$ a.s., this implies the second part of the lemma.  □

The following lemma is needed to verify condition (C4).

LEMMA A.14.  *Assume that conditions* (D1) *and* (D2) *hold and* $a_q^2 + \sigma_q^2 > 0$. *Then*

$$\Phi(x, t) = \tilde{\mathbf{E}}_x f(x'\Pi_{2q+1}, t), \qquad x \in S, t \in \mathbb{R},$$

*is uniformly continuous on* $S$ *for every measurable bounded function* $f : S \times \mathbb{R} \to \mathbb{R}$; *that is,*

$$\lim_{\varepsilon \to 0} \sup_{|x-y| \leq \varepsilon} \sup_{t \in \mathbb{R}} |\Phi(x, t) - \Phi(y, t)| = 0.$$

PROOF. Let $V : \mathbb{R}^q \to [0, \infty)$ be a continuous function such that $V(z) = 0$ for $|z| \geq 1$ and $\int_{\mathbb{R}^q} V(z) \, dz = 1$. For some $\epsilon \in (0, 1)$, define $K_\epsilon = \{y \in \mathbb{R}^q : |\langle y \rangle_1| \geq \epsilon, |y| \leq 1/\epsilon\}$, $\nu_\epsilon = \epsilon/4$ and $g_\epsilon(x) = \int_{|y| \leq 1} \chi_{K_\epsilon}(x + \nu_\epsilon y) V(y) \, dy$. Then $g_\epsilon : \mathbb{R}^q \to [0, 1]$ is continuous, such that $g_\epsilon(x) \leq \chi_{K_{\epsilon/4}}(x)$ and $\overline{g}_\epsilon(x) = 1 - g_\epsilon(x) \leq \chi_{K_{4\epsilon}^c}(x)$ for every $x \in \mathbb{R}^q$. We can represent the function $\Phi$ in the following form:

$$\Phi(x, t) = \tilde{\mathbf{E}}_x f(x' \Pi_{2q+1}, t) = \tilde{\mathbf{E}}_x g_\epsilon(\rho(x)) f(x' \Pi_{2q+1}, t) + \Delta_\epsilon(x),$$

where $\Delta_\epsilon(x) = \tilde{\mathbf{E}}_x \overline{g}_\epsilon(\rho(x)) f(x' \Pi_{2q+1}, t)$. By (A.11), setting $f^* = \sup |f|$, we obtain

$$\Delta_\epsilon^* = \sup_{x \in S} |\Delta_\epsilon(x)| \leq f^* \sup_{x \in S} \tilde{\mathbf{P}}_x(\rho(x) \in K_{4\epsilon}^c) \to 0, \qquad \epsilon \to 0.$$

From the definition of $\tilde{\mathbf{E}}$ in (4.3) we obtain

$$\tilde{\mathbf{E}}_x g_\epsilon(\rho(x)) f(x' \Pi_{2q+1}, t) = \frac{1}{h(x)} \mathbf{E} g_\epsilon(\rho(x)) f_1(x' \Pi_{2q+1}, t),$$

where $f_1(z, t) = |z|^\lambda h(\overline{z}) f(z, t)$. By Lemma A.10 we can represent this term as

$$\mathbf{E} g_\epsilon(\rho(x)) f_1(x' \Pi_{2q+1}, t) = \mathbf{E} \int_{\mathbb{R}^q} \overline{f}_1(z, t) \psi_\epsilon(z, \rho(x)) \, dz$$

$$= \mathbf{E} \Psi_\epsilon(\rho(x), t)$$

with $\overline{f}_1(z, t) = \mathbf{E} f_1(z'T, t)$ and $\psi_\epsilon(z, \rho) = p_0(z, \rho) g_\epsilon(\rho)$. Here $\Psi_\epsilon$ allows the representation

$$\begin{aligned}
\Psi_\epsilon(\rho, t) &= \int_{\Omega_\delta} \overline{f}_1(z, t) \psi_\epsilon(z, \rho) \, dz + \int_{\Omega_\delta^c} \overline{f}_1(z, t) \psi_\epsilon(z, \rho) \, dz \\
(\text{A.12}) \qquad &= \Psi_{\epsilon, \delta}(\rho, t) + \Delta_{\epsilon, \delta}(\rho, t),
\end{aligned}$$

where $\Omega_\delta = \{y \in \mathbb{R}^q : \delta \leq |\langle y \rangle_j| \leq \delta^{-1}, j = 1, \ldots, q\}$. Next we show that for every $\epsilon > 0$,

$$(\text{A.13}) \qquad \lim_{\delta \to 0} \sup_{\rho \in K_{\epsilon/4}} \mathbf{P}(\theta(\rho) \in \Omega_\delta^c) = 0.$$

To this end note

$$\sup_{\rho \in K_{\epsilon/4}} \mathbf{P}(\theta(\rho) \in \Omega_\delta^c)$$

$$\leq \sum_{j=1}^{q} \sup_{|\langle \rho \rangle_1| \geq \epsilon/4} \mathbf{P}(|\theta_j(\rho)| < \delta) + \sup_{|\rho| \leq 4/\epsilon} \mathbf{P}(|\theta(\rho)| > 1/\delta)$$

$$\leq \sum_{j=1}^{q} \sup_{|\langle \rho \rangle_1| \geq \epsilon/4} \mathbf{P}(|\theta_j(\rho)| < \delta) + \delta \sup_{|\rho| \leq 4/\epsilon} \mathbf{E}|\theta(\rho)|.$$

By the definition of $\theta(\rho)$ in (A.1) we find for every $m > 0$ some constant $c_m > 0$ such that

$$\sup_{|\rho| \leq 4/\epsilon} \mathbf{E}|\theta(\rho)|^m \leq c_m/\epsilon^m < \infty.$$

Therefore the limit relation (A.2) implies (A.13). Moreover, notice that the last inequality yields

$$\lim_{N \to \infty} \sup_{|\rho| \leq 4/\epsilon} \mathbf{E}\chi_{\{|\theta(\rho)| > N\}}|\theta(\rho)|^\lambda = 0.$$

Next we estimate $\Delta_{\epsilon,\delta}(\rho,t)$ as defined in (A.12). Taking into account that

$$|\overline{f}_1(z,t)| \leq f^* h^* \mathbf{E}|T|^\lambda |z|^\lambda = f_1^* |z|^\lambda,$$

we obtain for $\rho \in \mathbb{R}^q$ and $N > 0$,

$$|\Delta_{\epsilon,\delta}(\rho,t)| \leq f_1^* g_\epsilon(\rho) \int_{\Omega_\delta^c} |z|^\lambda p_0(z,\rho)\,dz$$

$$= f_1^* g_\epsilon(\rho)\mathbf{E}|\theta(\rho)|^\lambda \chi_{\{\theta(\rho) \in \Omega_\delta^c\}}$$

$$\leq f_1^* \chi_{\{\rho \in K_{\epsilon/4}\}}(N^\lambda \mathbf{P}(\theta(\rho) \in \Omega_\delta^c) + \mathbf{E}\chi_{\{|\theta(\rho)| > N\}}|\theta(\rho)|^\lambda).$$

This together with (A.13) ensures for every $\epsilon > 0$,

$$\Delta_{\epsilon,\delta}^* = \sup_{\rho \in \mathbb{R}^q, t \in \mathbb{R}} |\Delta_{\epsilon,\delta}(\rho,t)| \to 0 \qquad \text{as } \delta \to 0.$$

From this we conclude for $x, y \in S$ such that $|x - y| \leq \eta$ and for $\mu > 0$,

$$|\Phi(x,t) - \Phi(y,t)|$$

$$\leq \mathbf{E}|h^{-1}(x)\Psi_{\epsilon,\delta}(\rho(x),t) - h^{-1}(y)\Psi_{\epsilon,\delta}(\rho(y),t)| + 2\Delta_\epsilon^* + 2\Delta_{\epsilon,\delta}^*$$

$$\leq \Psi_{\epsilon,\delta}^* |h^{-1}(x) - h^{-1}(y)| + 2\Delta_\epsilon^* + 2\Delta_{\epsilon,\delta}^*$$

$$+ \frac{q}{\delta^\lambda} f_1^* \text{mes}(\Omega_\delta)\mathbf{E}\sup_{z \in \Omega_\delta} |\psi_\epsilon(z,\rho(x)) - \psi_\epsilon(z,\rho(y))|,$$

where $\Psi_{\epsilon,\delta}^* = \sup |\Psi_{\epsilon,\delta}|$. Since $\psi_\epsilon(\cdot,\cdot)$ and $\rho(\cdot)$ are uniformly continuous on $\Omega_\delta \times \mathbb{R}^q$ and on $S$, respectively, taking the limits $\lim_{\epsilon \to 0} \lim_{\delta \to 0} \lim_{\eta \to 0}$ implies Lemma A.14. $\square$

**A.4. General Markov properties of $(x_n)_{n \in \mathbb{N}}$.**   We consider now the Markov chain $(x_n)_{n \in \mathbb{N}}$ as defined in (1.10). Criteria for uniform ergodicity are based on "small" sets. A set $\Gamma \in \mathcal{B}(S)$ is called a *small set* if there exists an $m \in \mathbb{N}$ and a nontrivial measure $\nu_m$ on $\mathcal{B}(S)$ [i.e., $\nu_m(S) > 0$] such that $\mathbf{P}^m(x,A) \geq \nu_m(A)$ for all $x \in \Gamma$ and $A \in \mathcal{B}(S)$. As a general reference on Markov processes, we refer to [18].

Lemma A.8. *Assume that conditions* (D1) *and* (D2) *hold,* $q \geq 2$ *and* $a_q^2 + \sigma_q^2 > 0$. *Then the following hold:*

(a) *The distribution of the random vector* $x_{2q+1}$ *has the following properties: let* $A$ *be a measurable set in* $S$ *and denote by* $\Lambda(\cdot)$ *the Lebesgue measure on* $\mathcal{B}(S)$*, then:*

(i) *if* $\Lambda(A) > 0$*, then* $\inf_{y \in S} \mathbf{P}_y(x_{2q+1} \in A) > 0$ *and* $\inf_{y \in S} \tilde{\mathbf{P}}_y(x_{2q+1} \in A) > 0$;

(ii) *if* $\Lambda(A) = 0$*, then* $\mathbf{P}_y(x_{2q+1} \in A) = 0$ *and* $\tilde{\mathbf{P}}_y(x_{2q+1} \in A) = 0$ *for all* $y \in S$.

(b) *The Markov chain* $(x_n)_{n \in \mathbb{N}}$ *(with respect to both measures* $\mathbf{P}$ *and* $\tilde{\mathbf{P}}$*) is* $\Lambda$*-irreducible and aperiodic. Moreover, every measurable subset of* $S$ *is small.*

Proof. (a) Recall that $x_n' = x'\Pi_n/|x'\Pi_n|$. Note that for every $x \in S$ and every measurable set $A \in S$,

$$\mathbf{P}_x(x_{2q+1} \in A) = \mathbf{P}(x'\Pi_{2q+1} \in B_A),$$

$$\tilde{\mathbf{P}}_x(x_{2q+1} \in A) = \tilde{\mathbf{P}}_x(x'\Pi_{2q+1} \in B_A),$$

where $B_A = L^{-1}(A) = \{y \in \mathbb{R}^q \setminus \{0\} : L(y) \in A\}$ and $L(y) = y/|y|$. From (A.9) we obtain for some $0 < \delta < 1$,

$$\text{(A.14)} \qquad \begin{aligned} \mathbf{P}^{2q+1}(x, A) &\geq p_*(\delta)\mu_\delta(B_A) = \nu_\delta(A), \\ \tilde{\mathbf{P}}^{2q+1}(x, A) &\geq \tilde{p}_*(\delta)\tilde{\mu}_\delta(B_A) = \tilde{\nu}_\delta(A) \end{aligned}$$

for positive constants $p_*(\delta)$ and $\tilde{p}_*(\delta)$.

Next we show

$$\text{(A.15)} \qquad \Lambda(A) > 0 \quad \Longrightarrow \quad \text{mes}(B_A) > 0.$$

Recall that $q \geq 2$, hence, if $\Lambda(A) > 0$, there exists a open set $V \subseteq A \subseteq S$ with $\Lambda(V) > 0$. Then $L^{-1}(V) \subseteq B_A$, but this set is open and nonempty in $\mathbb{R}^q$ [$L(\cdot)$ is a continuous function on $\mathbb{R}^q \setminus \{0\}$ and $V \subset L^{-1}(V)$], therefore $\text{mes}(L^{-1}(V)) > 0$, which gives (A.15). If $\text{mes}(B_A) > 0$, then, by Lemma A.13, there exists some $\delta > 0$ such that $\mu_\delta(B_A) > 0$ and $\tilde{\mu}_\delta(B_A) > 0$. Then (i) follows from (A.14). Next we show that

$$\text{(A.16)} \qquad \Lambda(A) = 0 \quad \Longrightarrow \quad \text{mes}(B_A) = 0.$$

Assume that $\text{mes}(B_A) > 0$. Then there exists an open set $V \subset B_A$ with $\text{mes}(V) > 0$. By definition of $B_A$ the image $U = L(V) = \{L(y)\ y \in V\} \subseteq A$. We show that $U$ is an open set in $S$. Indeed, for $z_0 \in U$ there exists $y_0 \in V$ such that $z_0 = L(y_0) = y_0/|y_0|$. Since $V$ is open, there exists some $\delta > 0$ such that $\{y \in \mathbb{R}^q : |y - y_0| < \delta\} \subset V$. Set $\varepsilon = \delta/|y_0|$ and take $z \in S$ such that $|z - z_0| < \varepsilon$. Note that for $y_z = |y_0|z$ we have $L(y_z) = z$ and

$$|y_z - y_0| = |y_0||z - z_0| < |y_0|\varepsilon = \delta.$$

Hence, $y_z \in V$ and therefore $z \in U$, that is, $\{z \in S : |z - z_0| < \varepsilon\} \subset U$. Consequently $U = L(V)$ is an open set in $S$. For $q \geq 2$, the Lebesgue measure of any open nonempty set in $S$ is positive. This is a contradiction to $\Lambda(A) = 0$ and, hence, (A.14) holds. Furthermore, if $\mathrm{mes}(B_A) = 0$, then by Lemma A.10 and Corollary A.11,

$$\mathbf{P}_y(x_{2q+1} \in A) = \mathbf{E}\mathbf{P}(y'\Pi_{2q+1} \in B_A | \rho(y))$$
$$= \mathbf{E}\int_{B_A} p_1(z|\rho(y)) \, dz = 0,$$
$$\tilde{\mathbf{P}}_y(x_{q+1} \in A) = \tilde{\mathbf{E}}_y\tilde{\mathbf{P}}_y(y'\Pi_{2q+1} \in B_A | \rho(y))$$
$$= \tilde{\mathbf{E}}_y\int_{B_A} \tilde{p}_1(z|\rho(y)) \, dz = 0.$$

(b) Note that (i) and (ii) immediately imply $\Lambda$-irreducibility and aperiodicity. From inequalities (A.14) we conclude then that every measurable subset in $S$ is small. $\square$

LEMMA A.9.   *Assume that conditions* (D1) *and* (D2) *hold,* $q \geq 2$ *and* $a_q^2 + \sigma_q^2 > 0$. *Then the Markov chain* $(x_n)_{n \geq 0}$ *with state space* $S$ *is positive Harris recurrent and uniformly geometric ergodic with respect to* $\mathbf{P}$ *(and* $\tilde{\mathbf{P}}$). *It has invariant measure* $\pi(\cdot)$ *[and* $\tilde{\pi}(\cdot)$, *resp.], which is equivalent to Lebesgue measure* $\Lambda(\cdot)$ *on* $S$.

PROOF.   Define $V : \mathbb{R}^q \to [1, \infty)$ by $V(y) = 1 + |\langle y \rangle_1|$. Then

$$\mathbf{E}_xV(x_1) = 1 + \mathbf{E}\varsigma(x) = L(x)V(x),$$

where $\varsigma(x) = |\langle x'A_1 \rangle_1|/|x'A_1|$ and $L(x) = (1 + \mathbf{E}\varsigma(x))/V(x)$. Since $a_q^2 + \sigma_q^2 > 0$ implies that $\alpha_{q1}^2 > 0$, $\mathbf{P}$-a.s., we obtain

$$\lim_{|\langle x \rangle_1| \to 1 : x \in S} L(x) = \frac{1}{2}\left(1 + \mathbf{E}\frac{|\alpha_{11}|}{|\alpha_1|}\right)$$
$$\leq \frac{1}{2}\left(1 + \mathbf{E}\frac{|\alpha_{11}|}{\sqrt{\alpha_{11}^2 + \alpha_{q1}^2}}\right) < 1.$$

Thus, there exist $r > 0$ and $\epsilon < 1$ such that $\sup_{|\langle x \rangle_1| > r} L(x) < 1 - \epsilon$, and we obtain that $V(\cdot)$ satisfies on the set $\Gamma = \{x \in S : |\langle x \rangle_1| \leq r\}$:

$$\sup_{x \in \Gamma} \int_S V(y)p(x, dy) < \infty$$

and, for some $\epsilon \in (0,1)$,

$$\int_S V(y)p(x,dy) < (1-\epsilon)V(x) \qquad \text{for all } x \in \Gamma^c.$$

By the second part of Lemma A.8 every subset of $S$ is small. Since $(x_n)_{n\geq 0}$ is aperiodic, $(x_n)_{n\geq 0}$ is uniformly geometric ergodic with respect to $\mathbf{P}$ (see [18], page 355). In the same way uniform geometric ergodicity of $(x_n)_{n\geq 0}$ with respect to $\tilde{\mathbf{P}}$ can be shown. Therefore, $(x_n)_{n\geq 0}$ has stationary distributions $\pi(\cdot)$ and $\tilde{\pi}(\cdot)$, respectively. Next we use Lemma A.8(a) to show that $\pi$, respectively, $\tilde{\pi}$ are equivalent to Lebesgue measure on $S$. If $\pi(A) = \lim_{n\to\infty} \mathbf{P}_x(x_n \in A) = 0$ and $\Lambda(A) > 0$, then by Lemma A.8(a)(i), we obtain the following contradiction

$$\pi(A) = \lim_{n\to\infty} \mathbf{P}_x(x_{n+2q+1} \in A) = \lim_{n\to\infty} \int_S \mathbf{P}_y(x_{2q+1} \in A)\mathbf{P}^{(n)}(x,dy)$$

$$\geq \inf_{y\in S} \mathbf{P}_y(x_{2q+1} \in A) > 0.$$

Next, if $\Lambda(A) = 0$, then by Lemma A.8(a)(ii),

$$\pi(A) = \lim_{n\to\infty} \mathbf{P}_x(x_{n+2q+1} \in A)$$

$$= \lim_{n\to\infty} \int_S \mathbf{P}_y(x_{2q+1} \in A)\mathbf{P}^{(n)}(x,dy) = 0.$$

Hence, $\pi(\cdot)$ and $\Lambda(\cdot)$ are equivalent on $S$. In the same way we obtain the equivalence of $\tilde{\pi}(\cdot)$ and $\Lambda(\cdot)$ on $S$. □

### A.5. A property of $\psi$.

LEMMA A.10. *If conditions* (D0) *and* (D4) *hold, then the function* $\psi(x,u)$ *defined in* (4.11) *is nonnegative, and for all* $x = (\langle x\rangle_1, \ldots, \langle x\rangle_q)' \in S$ *with* $\langle x\rangle_1 \neq 0$,

$$(A.17) \qquad \text{mes}(\{u \geq 0 : \psi(x,u) > 0\}) > 0,$$

*where* $\text{mes}(\cdot)$ *denotes Lebesgues measure on* $\mathbb{R}$.

PROOF. By definition we have $\psi(x,u) = \mathbf{P}(\tau_1 + \tau_2 > u) - \mathbf{P}(\tau_1 > u)$ with $\tau_1 = x'A_1Y_1$ and $\tau_2 = x'\zeta_1 = \langle x\rangle_1 \xi_1$. If $\langle x\rangle_1 = 0$, then $\tau_2 = 0$, and therefore $\psi_0(x,u) = 0$. We show that $\psi_0(x,u) \geq 0$ if $\langle x\rangle_1 \neq 0$. By conditioning on $\tau_2$ we get

$$\psi(x,u) = \int_0^\infty (\mathbf{P}(u - t < \tau_1 \leq u) - \mathbf{P}(u < \tau_1 \leq u + t))p_{\tau_2}(t)\,dt$$

$$= \int_0^\infty \delta(u,t)p_{\tau_2}(t)\,dt,$$

where $p_{\tau_2}(\cdot)$ is the density of $\tau_2$, which is by condition (D4) symmetric and nonincreasing on $[0, \infty)$. Setting $\mathcal{A} = \sigma\{A_i, i \in \mathbb{N}\}$, again by condition (D4), the conditional density $p_{\tau_1}(\cdot|\mathcal{A})$ of $\tau_1$ is symmetric and nonincreasing on $\mathbb{R}_+$ a.s. Therefore the nonconditional density $p_{\tau_1}(\cdot)$ of $\tau_1$ have the same properties. Thus for $0 \leq t \leq u$, we have

$$\delta(u, t) = \int_{u-t}^{u} p_{\tau_1}(a)\, da - \int_{u}^{u+t} p_{\tau_1}(a)\, da$$

$$= \int_{u-t}^{u} (p_{\tau_1}(a) - p_{\tau_1}(a+t))\, da \geq 0.$$

On the other hand, for $t > u$, we get

$$\delta(u, t) = \int_{u-t}^{0} p_{\tau_1}(a)\, da + \int_{0}^{u} p_{\tau_1}(a)\, da - \int_{u}^{u+t} p_{\tau_1}(a|\mathcal{A})\, da$$

$$= \int_{0}^{t-u} (p_{\tau_1}(a) - p_{\tau_1}(a+2u))\, da$$

$$\quad + \int_{0}^{u} (p_{\tau_1}(a) - p_{\tau_1}(a+u))\, da$$

$$\geq 0,$$

again since $p_{\tau_1}(\cdot|\mathcal{A})$ is nonincreasing on $\mathbb{R}_+$. This proves the first part of the lemma.

We show now (A.17). Let $a_0 > 0$ such that $p_{\tau_1}(a_0 - s) > p_{\tau_1}(a_0 + s)$ for every $0 < s < a_0$ and $0 < t_0 < a_0$ such that $\mathbf{P}(\tau_2 > t_0) > 0$. Then for $t_0 < t < a_0$ and $a_0 < u < a_0 + t_0/2$,

$$\delta(u, t) = \int_{u-t}^{u} (p_{\tau_1}(a) - p_{\tau_1}(a+t))\, da$$

$$\geq \int_{a_0 - t_0/2}^{a_0} (p_{\tau_1}(a) - p_{\tau_1}(a+t_0))\, da$$

$$> 0.$$

This implies (A.17) immediately.  □

## REFERENCES

[1] BASRAK, B., DAVIS, R. D. and MIKOSCH, T. (2002). Regular variation of GARCH processes. *Stochastic Process. Appl.* **99** 95–115. MR1894253

[2] BINGHAM, N. H., GOLDIE, C. M. and TEUGELS, J. L. (1989). *Regular Variation*, revised paperback ed. Cambridge Univ. Press. MR1015093

[3] BORKOVEC, M. and KLÜPPELBERG, C. (2001). The tail of the stationary distribution of an autoregressive process with ARCH(1) errors. *Ann. Appl. Probab.* **11** 1220–1241. MR1878296

[4] BRANDT, A., FRANKEN, P. and LISEK, B. (1990). *Stationary Stochastic Models.* Wiley, Chichester. MR1086872

[5] DIACONIS, P. and FREEDMAN, D. (1999). Iterated random functions. *SIAM Review* **41** 45–76. MR1669737

[6] DIEBOLT, J. and GUEGAN, D. (1993). Tail behaviour of the stationary density of general non-linear autoregressive processes of order 1. *J. Appl. Probab.* **30** 315–329. MR1212664

[7] DUNFORD, N. and SCHWARTZ, J. T. (1958). *Linear Operators. Part I*: *General Theory*. Interscience Publishers, Inc., New York. MR117523

[8] ENGLE, R. F. (ed.) (1995). *ARCH. Selected Readings.* Oxford Univ. Press.

[9] FEIGIN, P. D. and TWEEDIE, R. D. (1985). Random coefficient autoregressive processes: A Markov chain analysis of stationarity and finiteness of moments. *J. Time Ser. Anal.* **6** 1–14. MR792428

[10] GOLDIE, C. M. (1991). Implicit renewal theory and tails of solutions of random equations. *Ann. Appl. Probab.* **1** 126–166. MR1097468

[11] GOLDIE, C. M. and MALLER, R. (2000). Stability of perpetuities. *Ann. Probab.* **28** 1195–1218. MR1797309

[12] HAAN, L., DE RESNICK, S. I., ROOTZÉN, H. and DE VRIES, C. G. (1989). Extremal behaviour of solutions to a stochastic difference equation, with applications to ARCH processes. *Stochastic Process. Appl.* **32** 213–224. MR1014450

[13] IBRAGIMOV, I. A. and HASMINSKII, R. Z. (1979). *Statistical Estimation*: *Asymptotic Theory*. Springer, Berlin. MR620321

[14] KESTEN, H. (1973). Random difference equations and renewal theory for products of random matrixes. *Acta Math.* **131** 207–248. MR440724

[15] KESTEN, H. (1974). Renewal theory for functional of a Markov chain with general state space. *Ann. Probab.* **2** 355–386. MR365740

[16] KLÜPPELBERG, C. and PERGAMENCHTCHIKOV, S. (2003). Renewal theory for functionals of a Markov chain with compact state space. *Ann. Probab.* **31** 2270–2300. MR2016619

[17] LE PAGE, E. (1983). Théorèmes de renouvellement pour les produits de matrixs aléatoires. Equations aux différences aléatoires. Publ. Sém. Math., Univ. Rennes. MR863321

[18] MEYN, S. and TWEEDIE, R. (1993). *Markov Chains and Stochastic Stability.* Springer, New York. MR1287609

[19] MIKOSCH, T. and STARICA, C. (2000). Limit theory for the sample autocorrelations and extremes of a GARCH(1, 1) process. *Ann. Statist.* **28** 1427–1451. (An extended version is available at www.math.ku.dk/~mikosch.) MR1805791

[20] NICHOLLS, D. F. and QUINN, B. G. (1982). *Random-Coefficient Autoregressive Models*: *An Introduction. Lecture Notes in Statist.* **11**. Springer, New York. MR671255

CENTER FOR MATHEMATICAL SCIENCES
MUNICH UNIVERSITY OF TECHNOLOGY
D-85747 GARCHING
GERMANY
E-MAIL: cklu@ma.tum.de
URL: www.ma.tum.de/stat/

LABORATOIRE DE MATHÉMATIQUES
    RAPHAËL SALEM
UMR CNRS 6085
SITE COLBERT
UNIVERSITÉ DE ROUEN
F-76821 MONT-SAINT-AIGNAN, CEDEX
FRANCE
E-MAIL: Serge.Pergamenchtchikov@univ-rouen.fr